\documentclass{article}

\usepackage{latex_style}

\begin{document}
\title{Whitham Deformations and the Space of Harmonic Tori in $\S^3$}
\author{Emma Carberry and Ross Ogilvie}
\date{}
\maketitle

\begin{abstract}
In this paper we investigate the space of harmonic maps from a 2-torus to $\S^3$ using the spectral curve correspondence and Whitham deformations. In an open and dense subset of a parameter space we find that the space of harmonic maps is smooth and has dimension two. We also show that the points that correspond to minimal tori (conformal harmonic maps) are either smooth points of dimension two or singular.
\end{abstract}

\section{Introduction}
In this paper we investigate the space of harmonic maps from a 2-torus to $\S^3$. From work of Hitchin\cite{Hitchin1990}, each harmonic map uniquely corresponds to spectral data $(\Sigma,\Theta^1,\Theta^2,E)$ consisting of a real hyperelliptic curve $\Sigma$, a pair of meromorphic differentials $\Theta^1, \Theta^2$, and a quaternionic line bundle $E$, which satisfy a litany of conditions (this correspondence excludes the case of conformal harmonic maps into a 2-sphere). In particular, there are constraints on the periods of the differentials. This correspondence makes it possible to construct a parameter space for the \emph{spectral triples} $(\Sigma, \Theta^1, \Theta^2)$.  Letting $\mathcal{M}_g$ be the space of spectral triples $(\Sigma,\Theta^1,\Theta^2)$ with a spectral curve $\Sigma$ of genus $g$,
the conditions that spectral data must satisfy allow us to parameterise spectral triples by polynomials of a fixed degree. We can then identify $\mathcal{M}_g$ as a subset of this parameter space. It is possible then to characterise the tangent vectors to $\mathcal{M}_g$ by equations~\eqref{eqn:EMPDi} and~\eqref{eqn:Q}, using Whitham deformations. At every point $p$ of an open and dense subset $\mathcal{U}_{ab}$ of the space of polynomials, where some additional assumptions hold, we find all solutions to these equations to demonstrate that $T_p\mathcal{M}_g$ is two-dimensional (Lemmata~\ref{lem:tangent generic} and~\ref{lem:tangent G}). Having established that the dimension is constant, it naturally follows that $\mathcal{M}_g\cap \mathcal{U}_{ab}$ is a manifold (Theorem~\ref{thm:moduli manifold}). A similar analysis at the conformal harmonic maps show that they are isolated in $\mathcal{M}_g$ or smooth points of dimension two.

The essential methodology is to describe the deformations that preserve the periods of differentials, the so called Whitham deformations. This method was first discovered for the Korteweg-de Vries equation~\cite{Flaschka1980,Lax1983}, before being developed generally for other integrable systems~\cite{Krichever1995,Gorsky1995}. 
It is difficult to describe the curves with a constraint on their periods directly, due to the Schottky problem~\cite{Grushevsky2012,Shiota1986}. Whitham deformations are therefore a powerful technique to extract local information.
The form of Whitham deformations used here resemble their application in the theory of constant mean curvature surfaces, with similar goals. In~\cite{Kilian2015a} it was shown that the space of equivariant CMC tori in $\S^3$ is a connected infinite graph.~\cite{Carberry2016a} 
uses Whitham deformations to show that for each $H > 0$ and each fixed genus of the spectral curve, spectral triples of tori of constant mean curvature $H$ in $\mathbb{S}^3 $ are dense amongst those of CMC planes. When the ambient space is $\mathbb{R}^3$, the tori are no longer dense and~\cite{Carberry2016} gives a partial description of their closure.

The reason to concentrate on spectral triples and ignore the line bundle is that deformations of the line bundle are well understood. Given a spectral datum, if one holds $\Sigma$ fixed then it is not possible to deform the differentials, because they must have integral periods, but it is possible to deform $E$. These are \emph{isospectral deformations} and the space of possible $E$ are described by a real $g$-dimensional torus in $\operatorname{Pic}^{g+1}(\Sigma)$. Thus we focus on \emph{non-isospectral deformations}, which are deformations of $\Sigma$ that preserve the integrality conditions of the differentials.

In our preliminaries (Section~\ref{sec:preliminaries}) we lay out the conditions that spectral data must satisfy and show how these conditions may be used to represent spectral triples $(\Sigma,\Theta^1,\Theta^2)$ as a triple of polynomials $(P,b^1,b^2)$. We describe an open and dense subspace $\mathcal{U}$ of the space of polynomials that partially satisfy the conditions to be spectral data. Hence $\mathcal{M}_g$ lies in $\mathcal{U}$. 
If we take a path in $\mathcal{M}_g$ parameterised by $t$, then because the periods of the differentials are integral, they are constant. Thus the $t$-derivative of each differential is exact, and yields a meromorphic function. Just as we characterised the differentials by a polynomial, so too can this pair of meromorphic functions be described by polynomials $(\hat{c}^1,\hat{c}^2)$.
The central result of Section~\ref{sec:deformations} is Lemma~\ref{lem:unique cs}, which shows each tangent vector to $\mathcal{M}_g$ is associated to a unique pair of polynomials $(\hat{c}^1,\hat{c}^2)$ and shows that this pair must satisfy a certain restrictive relationship~\eqref{eqn:Q reduced}.

A standalone interlude, Section~\ref{sec:bezout identity}, gives elementary but useful variations of B\'ezout's identity adapted for the types of polynomial equations found in this field of study. With these tools, it is then possible to establish the main result of the paper in Section~\ref{sec:tangent space}. This is the result that $\mathcal{M}_g\cap\mathcal{U}_{ab}$ is a smooth two-dimensional manifold, for $\mathcal{U}_{ab}$ an open and dense subset of the parameter space.
Section~\ref{sec:conformal} gives a similar result, that generically conformal harmonic maps (i.e. minimal surfaces) are smooth points of $M_g$.
% The final section contains a proof that the assumptions of Theorems~\ref{thm:moduli manifold} and~\ref{thm:conformal moduli manifold} hold automatically for spectral curves with low genus, and in higher genus discusses the obstacles to proving analogous results for the remaining points of $\mathcal{M}_g$.

%%%%%%%%%%%%%%%%%%%%%%%%%%%%%%%%%%%%%%%%%%%%%%%%%%%%%%%%%%%%%%%%%%%%
%%%%%%%%%%%%%%%%%%%%%%%%%%%%%%%%%%%%%%%%%%%%%%%%%%%%%%%%%%%%%%%%%%%%
%%%%%%%%%%%%%%%%%%%%%%%%%%%%%%%%%%%%%%%%%%%%%%%%%%%%%%%%%%%%%%%%%%%%
%%%%%%%%%%%%%%%%%%%%%%%%%%%%%%%%%%%%%%%%%%%%%%%%%%%%%%%%%%%%%%%%%%%%
\section{Preliminaries}\label{sec:preliminaries}
Hitchin~\cite{Hitchin1990} investigated harmonic maps from a torus to the $3$-sphere and characterised them in terms of a spectral curve construction. The spectral data corresponding to a harmonic map $T^2 \to \S^3$ consists of a tuple $(\Sigma, \Theta^1, \Theta^2, E)$. $\Sigma$ is called the \emph{spectral curve} and is a real hyperelliptic curve over $\CP^1$ with a pair of meromorphic differentials $\Theta^1,\Theta^2$ and line bundle $E$. Certain geometric features of a harmonic map may be discerned from properties of spectral data. For example, a harmonic map is conformal (and therefore minimal) exactly when the spectral curve is branched over $0$. Theorems of Hitchin~\cite[Theorems~8.1,~8.20]{Hitchin1990} provide a correspondence between harmonic maps from the torus to $\S^3$ and tuples $(\Sigma,\Theta^1,\Theta^2,E)$ satisfying properties~\ref{P:real curve}--\ref{P:quaternionic}. The harmonic map is uniquely determined by its spectral data, up to the action of $\SO(4)$ on $\S^3$. We shall now enumerate those properties.

Suppose that $\Sigma$ has genus $g$ and is described by $\eta^2 = P(\zeta)$ in the total space of $\pi: \mathcal{O}(g+1) \to \CP^1$ with $P(\zeta)$ a section of $\mathcal{O}(2g+2)$. This curve is hyperelliptic with involution $\sigma(\zeta,\eta) = (\zeta,-\eta)$. 
% A real structure is an anitholomorphic involution. We will use the real structure $\zeta \mapsto \cji{\zeta}$ on $\CP^1$. 
A spectral curve has the following properties.
\begin{enumerate}[label=(P.\arabic*)]
    \item\label{P:real curve} Real curve: $P(\zeta)$ is a real section of $\mathcal{O}(2g+2)$.
    \item\label{P:no real zeroes} No real zeroes: $P(\zeta)$ has no zeroes on the unit circle $\S^1\subset\CP^1$.
    \item\label{P:simple zeroes} Simple zeroes: $P(\zeta)$ has only simple zeroes over $\zeta = 0,\infty$.
\end{enumerate}
% Let the holomorphic involution $\rho$ be the restriction of $\rho_{g+1}$ on $\mathcal{O}(g+1)$ to the spectral curve. The fixed points of $\sigma$ are the roots of $P$, so~\ref{P:no real zeroes} ensures that $\rho\circ \sigma$ is fixed point free.

\textbf{Assumption:} In this paper we consider only nonsingular spectral curves. This means $P(\zeta)$ may only have simple roots.

Sections of $\mathcal{O}(k)$ can be identified with polynomials of degree at most $k$. 
We use the real structure on the line bundle $\mathcal{O}(k)$ given by
\[
\rho_k : (\zeta,\eta) \mapsto \left(\cji{\zeta}, \bar{\zeta}^{-k}\bar{\eta}\right).
\]
\begin{defn}\label{def:mathcal P}
Let $\mathcal{P}^k$ be the space of polynomials of degree at most $k$. We define the real polynomials $\mathcal{P}^k_\R$ to be
\[
\mathcal{P}^k_\R
= \Set{ q_0 + \cdots + q_k \zeta^k \in \mathcal{P}^k }{ q_i = \bar{q_{k-i}} \text{ for all } 0 \leq i \leq k}.
\]
$\mathcal{P}^k_\R$ is a real vector space of dimension $k+1$.
\end{defn}
The polynomials in $\mathcal{P}^k_\R$ are those that correspond to sections invariant under pullback by the real involution $\rho_k$. 
Every root of a real polynomial must either lie on the unit circle or come in a conjugate-inverse pair. Identifying $P$ with a polynomial in $\mathcal{P}^{2g+2}_\R$, we fix the following scaling of $P$,
\[
    P(\zeta) = \prod_{i=0}^{g} (\zeta - \alpha_i)(1 -\bar{\alpha}_i \zeta),
    \labelthis{eqn:def P}
\]
for $\alpha_0,\ldots,\alpha_g$ in the unit disc. Roots on the unit circle are excluded by~\ref{P:no real zeroes}.
A nice feature of this scaling is that it is well-behaved if one branch point is zero; the corresponding factor becomes $\zeta$. A hyperelliptic curve is determined by its branch points in $\CP^1$, uniquely up to automorphism of $\CP^1$. Since we have distinguished the points $0,1$, and $-1$ we have fixed the automorphism, and so every spectral curve corresponds to a unique polynomial of the form~\eqref{eqn:def P}.

Each differential $\Theta^1,\Theta^2$ must satisfy the following conditions.
\begin{enumerate}[resume*]
\item\label{P:poles} Poles: The differentials have residue-free double poles at $\pi^{-1}\{0,\infty\}$ but are otherwise holomorphic.
\item\label{P:symmetry} Symmetry: The differentials satisfy $\sigma^* \Theta^i = - \Theta^i$.
\item\label{P:reality} Reality: The differentials satisfy $\rho_{g+3}^* \Theta^i = - \bar{\Theta}^i$.
% \item\label{P:imaginary periods} Imaginary Periods: The differentials have purely imaginary periods.
\item\label{P:linear independence} Linear independence: The principal parts of the differentials $\Theta^1$ and $\Theta^2$ are linearly independent over $\R$.
\end{enumerate}

Every curve satisfying~\ref{P:real curve}--\ref{P:simple zeroes} has differentials that satisfy~\ref{P:poles}--\ref{P:reality}, necessarily of the form
\[
\Theta = b(\zeta)\frac{d\zeta}{\zeta^2\eta} = \left( b_0 + b_1 \zeta + \cdots + b_{g+3} \zeta^{g+3} \right) \frac{d\zeta}{\zeta^2\eta},
\labelthis{eqn:def b}
\]
for some $b(\zeta) \in \mathcal{P}^{g+3}_\R$, by~\cite[Prop~III.1.10]{Miranda1995}. If $\Sigma$ does not have a branch point at $0$, then the residue of $\Theta$ is
\[
\operatorname*{res}_{\zeta=0}\Theta = b_1 - \frac{1}{2}\frac{P_1}{P_0}b_0,
\labelthis{eqn:residue}
\]
where subscripts denote coefficients of the polynomials. This quantity must therefore vanish by~\ref{P:poles}. On the other hand, if $P_0$ is zero then the same condition forces $b_0=0$. In light of this, we may rephrase equation~\eqref{eqn:residue},
\[
P_1b_0 - 2P_0b_1 = 0.
\labelthis{eqn:residue condition}
\]

The differentials belonging to spectral data however have properties further than~\ref{P:poles}--\ref{P:linear independence}, which are `hard' to satisfy. 
\begin{enumerate}[resume*]
\item\label{P:periods} Periods: The periods of the differentials $\Theta^1$ and $\Theta^2$ lie in $2\pi\iu\Z$.
\item\label{P:closing} Closing conditions: 
% Let $\mu^i$ be a meromorphic function on $\Sigma\setminus \pi^{-1}\{0,\infty\}$ such that $\Theta^i = d\log \mu^i$ and $\mu^i\sigma^*\mu^i = 1$. Then $\mu^i$ has value $1$ at $\pi^{-1}\{\pm 1\}$
Let $\gamma_+$ be a path in $\Sigma$ connecting the two points of $\pi^{-1}(1) = \{\xi_1, \sigma(\xi_1)\}$ and likewise let $\gamma_-$ connect the two points $\{\xi_{-1}, \sigma(\xi_{-1})\}$ over $-1$. Then
\[
\int_{\gamma_+} \Theta^i \in 2\pi\iu\Z,\;\;\text{ and }
\int_{\gamma_-} \Theta^i \in 2\pi\iu\Z.
\]
\end{enumerate}
Finally, the line bundle of the spectral data has the following property.
\begin{enumerate}[resume*]
\item\label{P:quaternionic} Quaternionic: $E^*$ is a line bundle of degree $g+1$ that is quaternionic with respect to the involution $\rho_{g+1}\circ \sigma$.
\end{enumerate}

If the choice of curve and differentials is fixed, one is free to choose $E$ subject only to~\ref{P:quaternionic}. There are many such choices; they form a real $g$ dimensional torus in the Jacobian of $\Sigma$. Variations of $E$ alone are called isospectral deformations. Conversely, if we have a triple $(\Sigma,\Theta^1,\Theta^2)$ satisfying the above conditions, then there always exists such a line bundle $E$ completing the tuple. Thus we focus our attention on the problem of deforming the spectral triple $(\Sigma,\Theta^1,\Theta^2)$; so called non-isospectral deformations.

We have seen that any spectral triple $(\Sigma,\Theta^1,\Theta^2)$ may be described in terms of polynomials $(P,b^1,b^2) \in \mathcal{P}^{2g+2}_\R\times\mathcal{P}^{g+3}_\R\times\mathcal{P}^{g+3}_\R$.
\begin{defn}\label{def:M_g}
We call $\mathcal{M}_g$ the space of spectral triples $(\Sigma,\Theta^1,\Theta^2)$ satisfying conditions~\ref{P:real curve}--\ref{P:closing} and such that $\Sigma$ has genus $g$ and is nonsingular. 

Let $\mathcal{U}$ be the following open and dense subset of $\mathcal{P}^{2g+2}_\R\times\mathcal{P}^{g+3}_\R\times\mathcal{P}^{g+3}_\R$. If $(P,b^1,b^2)$ is a point of $\mathcal{U}$ then $P$ has no zeroes on the unit circle (cf.~\ref{P:no real zeroes}) and only simple zeroes elsewhere (cf.~\ref{P:simple zeroes} and $\Sigma$ nonsingular). The polynomials $b^i$ have at most a simple root at $\zeta=0$ (cf.~\ref{P:poles}).

We view $\mathcal{M}_g$ as a subset of $\mathcal{U}$ using the above correspondence between spectral triples and polynomials.
\end{defn}

\emph{A comment about notation.} 
% Throughout we will use diacritical marks to indicate the factors that a polynomial does or does not have.
A polynomial with a circumflex (hat) will be shown to have a factor of $\zeta^2-1$, and a tilde will indicate that common factors have been removed, cf.~\eqref{eqn:common factors}. 
We shall use a dash to denote differentiation with respect to $\zeta$ and a dot for differentiation with respect to $t$ evaluated at $t=0$.
When giving the solutions to equations we will use bold to signify a particular solution, which may or may not be unique.
% , whereas a solution without bold signifies any solution from the set of potential solutions. 
Given a tuple of polynomials, such as $(X,Y)$, we also give their degrees as a tuple, e.g. $(x,y)$ for $x=\deg X$ and $y=\deg Y$. 
Finally, we shall use $i$ and $j$ for indices ranging over $\{1,2\}$, with the understanding that they are not equal. 
% For example, if $i=1$, then we take $j=2$ and vice versa.

%%%%%%%%%%%%%%%%%%%%%%%%%%%%%%%%%%%%%%%%%%%%%%%%%%%%%%%%%%%%%%%%%%%%
%%%%%%%%%%%%%%%%%%%%%%%%%%%%%%%%%%%%%%%%%%%%%%%%%%%%%%%%%%%%%%%%%%%%
%%%%%%%%%%%%%%%%%%%%%%%%%%%%%%%%%%%%%%%%%%%%%%%%%%%%%%%%%%%%%%%%%%%%
%%%%%%%%%%%%%%%%%%%%%%%%%%%%%%%%%%%%%%%%%%%%%%%%%%%%%%%%%%%%%%%%%%%%
\section{Deformations of Harmonic Maps}\label{sec:deformations}

In this section we lay out the properties of infinitesimal deformations of spectral triples $(\Sigma,\Theta^1,\Theta^2)$ within the space $\mathcal{M}_g$. A deformation of spectral data is a path $\ell:(-\epsilon,\epsilon) \to \mathcal{M}_g$, $t \mapsto (P(t,\zeta),b^1(t,\zeta),b^2(t,\zeta))$. An infinitesimal deformation is the tangent vector of such a curve at $t=0$.
Denote $\Sigma = \Sigma(0)$ and likewise $b^i(\zeta) = b^i(0,\zeta)$.
More generally, omission of the parameter $t$ will correspond to evaluation at the point $t=0$.

Consider the $t$-derivative of $\Theta^i$. The periods of $\Theta^i$ are constant in $t$ by~\ref{P:periods}, so $\dot{\Theta}^i$ is exact. Write $\dot{\Theta^i} = d\theta^i$, for $\theta^i$ a meromorphic function on $\Sigma$. In order to write $\theta^i$ analogously to how we wrote $\Theta^i$ in~\eqref{eqn:def b} we need to know how a function may acquire additional poles when it is differentiated with respect to $t$. This must be handled delicately; consider the example of $f(t)(z) = (z^n+t)/{(z-t)}^n$, where $f(0) \equiv 1$ but $\dot{f}$ has a pole of order $n$ at $z=0$.

\begin{lem}\label{lem:pole order}
Consider a smooth family of hyperelliptic curves $\Sigma(t)$ and 
a smooth family of meromorphic functions $f(t)$, for $t\in(-\epsilon,\epsilon)$.
Fix a point $p\in\Sigma(0)$ that lies over $\beta \in \CP^1$.

Suppose $p$ is a not ramification point. If $f(t)$ is holomorphic for every $t$ then $\dot{f}$ is holomorphic. If ${(\zeta - \beta)}^n f(t)$ is holomorphic for every $t$ then $\dot{f}$ may have a pole of at most order $n$ at $p$. 

Suppose that $p$ is a ramification point. If $f(t)$ is holomorphic for every $t$ then $\dot{f}$ may have at worst a simple pole at $p$.
%  If ${(\zeta - \beta)}^n f(t)$ is holomorphic for every $t$ then $\dot{f}$ has a pole of at most order $2n+1$ at $p$.

\begin{proof}
Let us consider the unramified case first. Near $p$ we may write our function $f(t)(\zeta) = u(t,x,y) + \iu v(t,x,y)$ for $\zeta = x + \iu y$. As $f(t)$ is holomorphic we have the Cauchy-Riemann equations $\partial_x u = \partial_y v$ and $\partial_y u = -\iu \partial_x v$. Because $f(t)$ is a smooth function we may differentiate by $t$ and interchange the order, yielding $\partial_x\partial_t u = \partial_y\partial_t v$ and $\partial_y\partial_t u = -\iu \partial_x\partial_t v$. This shows that $\dot{f}$ is holomorphic too. The second part then follows because the left hand side of
\[
\frac{d}{d t}{(\zeta - \beta)}^n f = {(\zeta - \beta)}^n \partial_t f,
\]
is now known to be holomorphic.

We move now to the ramified case. Let $\alpha(t)$ be a branch point of $\Sigma(t)$ such that $\alpha(0)=\beta$. Note that there is a unique choice. On $\Sigma(t)$ take $\xi(t)^2 = \zeta - \alpha(t)$ as local coordinate, so we may compute
% . In particular, $\dot{\xi} = -\tfrac{1}{2}\dot{\alpha}\xi^{-1}$ and
\[
\frac{d}{dt}f(t)(\xi)
= \frac{\partial f}{\partial t} - \frac{1}{2}\dot{\alpha}\xi^{-1}\frac{\partial f}{\partial\xi}.
\]
If $f(t)$ was holomorphic, then we see that $\dot{f}$ may have a simple pole. 
% If $(\zeta-\beta)^n f(t)$ was holomorphic, then its $t$-derivative may have a simple pole. Hence $\dot{f}$ may have a pole of $2n+1$, since $(\zeta-\beta)^{-n}$ has a pole of order $2n$.
\end{proof}
\end{lem}

As $\Theta^i$ has double poles without residues over $\zeta=0$ and $\infty$, it follows that $\theta^i$ may have simple poles at the roots of $P$ in $\C^\times$. 
If the curve $\Sigma$ is unbranched over $\zeta=0$ and $\infty$, then  $\theta^i$ has at worst simple poles there.
If the curve $\Sigma$ is branched over $\zeta=0$, then by~\eqref{eqn:residue condition} there exists a root $\beta(t)$ of $\Theta^i(t)$ such that $\zeta^2(\zeta-\beta(t))^{-1}\Theta^i(t)$ is holomorphic near $\zeta=0$ for all $t$ and $\beta(0)=0$. Differentiating gives
\[
\dot{\Theta^i} = \dot{\beta}\zeta^{-1}\Theta^i 
+ \zeta^{-1}\left.\frac{d}{dt} \zeta^2(\zeta-\beta(t))^{-1}\Theta^i(t) \right|_{t=0},
\]
from which we can see that $\dot{\Theta}^i$ has at worst a fourth order pole over $\zeta=0$ in this case. Therefore $\theta^i$ may have a triple pole there. 
Hence $\zeta\eta\theta^i$ is holomorphic over $\C$ and
for some degree $g+3$ polynomial $\hat{c}^i$
\[
\theta^i = \frac{1}{\zeta\eta}\hat{c}^i(\zeta).
\labelthis{eqn:def q dot}
\]

By definition, $\dot{\Theta}^i = d\theta^i$. This provides equations linking $\dot{b}^i$ and $\hat{c}^i$, 
\[
\dot{P} b^i - 2P\dot{b}^i = 2P\left( \hat{c}^i - \zeta\hat{c}^{i\prime} \right) + P'\zeta\hat{c}^i. \labelthis{eqn:EMPDi}
\]

It is natural to ask to what extent the functions $\theta^i$ determine an infinitesimal deformation $(\dot{P},\dot{b}^1,\dot{b}^2)$. In Lemma~\ref{lem:unique cs} we will show that $\theta^1$ and $\theta^2$ are uniquely determined by an infinitesimal deformation. Thus we may use these functions to identify the space of infinitesimal deformations that preserve periods with a subset of $\mathcal{P}^{g+3}\times \mathcal{P}^{g+3}$. The converse question is more difficult and is addressed in Section~\ref{sec:tangent space}.

\begin{lem}\label{lem:unique cs}
Given a point $(P,b^1,b^2)$ of $\mathcal{M}_g$ for which there are deformations, the polynomials $(\hat{c}^1,\hat{c}^2)$ are determined uniquely by a tangent vector $(\dot P, \dot b^1, \dot b^2)$ to $\mathcal{M}_g$.

\begin{proof}
Since the equations~\eqref{eqn:EMPDi} are linear in the components of the tangent vector, we need only demonstrate that the zero tangent vector uniquely corresponds to $\hat{c}^i = 0$. For the zero tangent vector,
\[
0 = 2P\left( \hat{c}^i - \zeta\hat {c}^{i\prime}\right) + P' \zeta\hat {c}^i.
\labelthis{eqn:EMPDi zero vector}
\]
The polynomial $P$ is either of the form $L$ or $\zeta L$, where $L$ has only roots in $\C^\times$. $L$ has degree either $2g+2$ or $2g$ respectively.
The assumption of a nonsingular spectral curve requires that $P$ and $P'$ have no common factors, hence evaluation of~\eqref{eqn:EMPDi zero vector} at any root $\alpha$ of $L$ shows that $P'(\alpha)\alpha\hat{c}^i(\alpha) = 0$, and hence $\alpha$ is a root of $\hat{c}^i$. This shows that $L$ divides $\hat{c}^i$, and for $g\geq 4$ the inequality $\deg L \geq 2g \geq g+4$ is sufficient to show that $\hat{c}^i$ is the zero polynomial, as it is a degree $g+3$ polynomial that is divisible by a polynomial of greater degree. To handle the remaining cases, $g<4$, we substitute in this factorisation of $\hat{c}^i$ and then remove the factor of $L$,
\begin{align*}
    % 0
    % &= 2L\frac{P}{L}\bra{ L\frac{\hat{c}^i}{L} - \zeta\left[ L' \frac{\hat {c}^i}{L} + L \bra{\frac{\hat {c}^i}{L}}' \right] }
    % + \left[L' \frac{P}{L} + L\bra{\frac{P}{L}}' \right] \zeta L \frac{\hat{c}^i}{L}  \\
    0 &= L \left [ 2\frac{P}{L}\frac{\hat{c}^i}{L} - 2\zeta\frac{P}{L}\bra{\frac{\hat {c}^i}{L}}' + \zeta\bra{\frac{P}{L}}'\frac{\hat {c}^i}{L} \right]
    -\zeta L' \frac{P}{L} \frac{\hat {c}^i}{L}.
\end{align*}
Again, this shows that $L$ divides $\hat{c}^i/L$. If $\deg L = 2g+2$, this shows that $\hat{c}^i$ is divisible by a polynomial of degree $4g+4$, and so must be zero for any $g$. If $\deg L = 2g$, then we have only shown that $\hat{c}^i$ vanishes for $g \geq 2$. We treat the two remaining cases, $\deg L = 2g$ and $g = 0$ or $1$, individually.

If $\deg L = 2g$ and $g = 1$, then $\hat{c}^i$ is a scalar multiple of $L^2$. Let $\hat{c}^i = a L^2$ and equation~\eqref{eqn:EMPDi zero vector} simplifies to
\[
0 = 2a L^3 \bra{ 3L - 2\zeta L' },
\]
which forces $a = 0$. For the second case if $g=0$ then $P=\zeta$ and $\hat{c}^i$ is a cubic polynomial. After removing the factor of $\zeta$,
\[
0
= 2(\hat{c}^i - \zeta\hat{c}^{i\prime}) + \hat{c}^i
= 3\hat{c}^i_0 + \hat{c}^i_1 \zeta - \hat{c}^i_2 \zeta^2 - 3\hat{c}^i_3 \zeta^3,
\]
which again shows that $\hat{c}^i$ is zero. Hence, the polynomials $\hat{c}^i$ are uniquely determined by a tangent vector $(\dot P, \dot b^1, \dot b^2)$ to $\mathcal{M}_g$ as claimed.
\end{proof}
\end{lem}

Continuing our line of inquiry into the properties of $\hat{c}^i$, recall our supposition that we are at a point $(P,b^1,b^2)$ of $\mathcal{M}_g$ that admits a deformation, from which we have defined polynomials $\hat{c}^1$ and $\hat{c}^2$ and derived the pair of equations~\eqref{eqn:EMPDi}. Note that the two equations~\eqref{eqn:EMPDi} for $i=1,2$ are not independent of one another, for they both contain $P$ and its derivatives. If we multiply the equations by $\hat c^2$ and $\hat c^1$ respectively and take the difference, we observe
\[
\dot P (b^1\hat c^2 - b^2\hat c^1) =  2P(\dot b^1\hat c^2 - \dot b^2\hat c^1 - \zeta\hat{c}^{1\prime}\hat c^2 + \zeta\hat{c}^{2\prime}\hat c^1).\labelthis{eqn:diffed}
\]
We will prove that $b^1\hat c^2 - b^2\hat c^1$ is divisible by $P$ by showing that it vanishes at every root of $P$. If $\alpha$ is a root of $P$ and not a root of $\dot{P}$, we see it is a root of $b^1\hat c^2 - b^2 \hat c^1$ immediately from~\eqref{eqn:diffed}. Suppose then that $P$ and $\dot P$ have a common root $\alpha$. If $\alpha=0$, then we know from~\eqref{eqn:residue condition} that $b^i_0=0$ and so $\zeta$ divides $b^i$. If $\alpha\neq 0$, from~\eqref{eqn:EMPDi} we have that
\begin{align*}
\dot P(\alpha) b^i(\alpha) &= 2P(\alpha)\left( \dot b^i(\alpha) + \hat c^i(\alpha) - \alpha\hat{c}^{i\prime}(\alpha)\right) +P'(\alpha)\alpha\hat c^i(\alpha) \\
0 &= 0 + P'(\alpha)\alpha\hat{c}^i(\alpha)
\end{align*}
But the assumption that the spectral curve is nonsingular forces $P'(\alpha)\neq 0$. Thus we may conclude that $\hat{c}^i(\alpha)=0$. Hence $P$ divides $b^1\hat c^2 - b^2 \hat c^1$ and there is some polynomial $\hat{Q}$ of degree at most four such that
\[
b^1 \hat{c}^2 - b^2 \hat{c}^1 = \hat{Q} P.
\labelthis{eqn:Q hat}
\]

Thus far we have only placed two conditions on the points along the deformation. First, that it must preserve the integral periods of $\Theta^1$ and $\Theta^2$, which allowed us to produce well-defined meromorphic functions $\theta^i$. And second that the differentials must have double poles over $\zeta=0$ and $\infty$ (Property~\ref{P:poles}), which allowed us to write $\theta^i$ as the quotient of a polynomial $\hat{c}^i$ by $\zeta\eta$. 

Applying~\ref{P:reality} forces the polynomials $\hat{c}^i$ to be imaginary (that is, $\iu\hat{c}^i \in \mathcal{P}^{g+3}_\R$). Next, consider the closing condition~\ref{P:closing}. Differentiating
\[
0 
= \left.\frac{d}{dt}\right|_{t=0} \int_{\gamma_+} \Theta^i
% = \int_{\gamma_+}\dot{\Theta}^i
= \int_{\gamma_+} d\theta^i
= \theta^i(\sigma(\xi_1)) - \theta^i(\xi_1)
\]
where $\xi_1$ is a point in $\Sigma$ over $\zeta=1$.
But
\[
\theta^i(\sigma(\xi_1)) = \sigma^* \theta^i (\xi_1) = - \theta^i(\xi_1).
\]
Thus $\hat c^i$ has a root at $\zeta=1$. The same reasoning applied to $\gamma_-$ leads to a root at $\zeta=-1$. Therefore let $\hat c^i(\zeta) = (\zeta^2 - 1) c^i(\zeta)$ for some $c^i \in \mathcal{P}^{g+1}_\R$.

As $\zeta^2-1$ is a factor of both polynomials $\hat{c}^i$, and $P$ has no zeroes on the unit circle, it follows from~\eqref{eqn:Q hat} that $\zeta^2-1$ must be a factor of $\hat{Q}$. Define $\hat Q = (\zeta^2-1)Q$ to give
\[
b^1 c^2 - b^2 c^1 = Q P
\labelthis{eqn:Q}
\]
for some real quadratic polynomial $Q$. The importance of this equation is that it ensures that the solutions to the two equations~\eqref{eqn:EMPDi} are consistent with one another.

The final condition on the spectral data that we are yet to satisfy is condition~\ref{P:poles}: that the differentials $\Theta^1$ and $\Theta^2$ are residue free. We shall require $P_1(t)b^i_0(t) - 2P_0(t)b^i_1(t) = 0$ to hold at every point of the deformation (from~\eqref{eqn:residue condition}). Taking derivatives,
\[
\dot{P}_1 b^i_0 + P_1 \dot{b}^i_0 - 2 \dot{P}_0 b^i_1 - 2 P_0 \dot{b}^i_1 = 0 \labelthis{eqn:residueTangent}
\]
holds for any for tangent vector $(\dot{P}, \dot{b}^1, \dot{b}^2)$ to $\mathcal{M}_g$.

\begin{SCfigure}
\begin{tikzpicture}
\clip (-4,-4) rectangle (4,4);

\draw (-1,-1) circle [radius=2];
\draw (1,-1) circle [radius=2];
\draw (0,0.73) circle [radius=2];

\node at (0,-0.135) {$F$};
\node at (-1.2,0.2) {$F^1$};
\node at (1.2,0.2) {$F^2$};
\node at (0,-1.7) {$G$};

\node at (0,3) {$P$};
\node at (-2.8,-2.8) {$b^1$};
\node at (2.8,-2.8) {$b^2$};
\end{tikzpicture}
\caption{\label{fig:common factors}
The common factors of $P$, $b^1$ and $b^2$ are represented as overlaps between the three circles.}
\end{SCfigure}
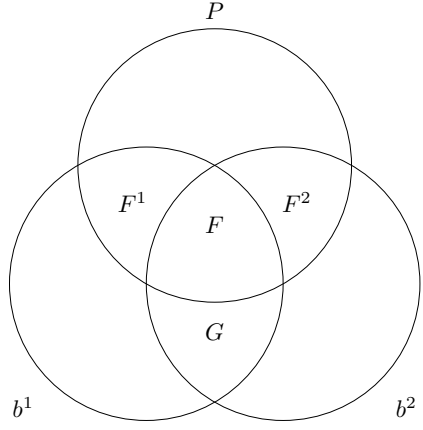

This covers the necessary properties that are shared by all infinitesimal deformations. However, there are further constraints in the special cases where the differentials have common zeroes or zeros at ramification points. In these cases, the polynomial $Q$ also shares a common root. Suppose we have the following common factors at the point $(P(0,\zeta),b^1(0,\zeta),b^2(0,\zeta)) \in \mathcal{M}_g$:
\begin{align*}
\gcd(P,b^1&, b^2) = F, \\
\labelthis{eqn:common factors}
\gcd(P/F, b^1/F) = F^1&,\;\; \gcd(P/F,b^2/F) = F^2,\\
\gcd(b^1/FF^1&, b^2/FF^2) = G,
\end{align*}
where we first find the common factor of all three polynomials, then remove any further factors that the differentials and $P$ share, and then finally remove any remaining factors common to $b^1$ and $b^2$. An graphic representation of this process is given in Figure~\ref{fig:common factors}. We write
\[
P = F F^1 F^2 \tilde{P},\;\; b^1 = F F^1 G \tilde{b}^1,\;\; b^2 = F F^2 G \tilde{b}^2.
\]
Because the spectral curve is nonsingular, $P$ has no repeated factors, and so the polynomials $F$, $F^1$, $F^2$, $\tilde{P}$, $\tilde{b}^1$ and $\tilde{b^2}$ are pairwise coprime. The polynomials $b^1$ and $b^2$ may have higher order roots, so it is not possible to say if $G$ is coprime to $F$, $F^1$ or $F^2$. The common factor of $b^1$ and $b^2$, and therefore any differential on $\Sigma$ satisfying conditions~\ref{P:poles}--\ref{P:periods}, is $FG$. We denote the degrees of the polynomials $F$, $F^1$, $F^2$ and $G$ as $d_F$, $d_1$, $d_2$ and $d_G$ respectively.

Inserting these factorisations into~\eqref{eqn:EMPDi}, we observe that
\[
\dot{P} F F^i G \tilde{b}^i = 2 F F^1 F^2 \tilde{P} (\dot{b}^i + \hat{c}^i - \zeta\hat{c}^{i\prime}) + \zeta(\zeta^2-1 ) P'c^i.
\labelthis{eqn:EMPDi all factors}
\]
Again, by the assumption of that the spectral curve is nonsingular, $P'$ does not share any common factors with $P$. Further $P$ has no roots on the unit circle. Hence we see that $FF^i$ divides $\zeta c^i$. 
% Conversely, given an arbitrary $c^i$, it would not be possible to solve this equation for $(\dot{P},\dot{b}^i)$ unless $FF^i$ divides $\zeta c^i$, otherwise we would have a contradiction. 
% We would like to say that $FF^i$ divides $c^i$ alone, but because $F$ may have a factor of $\zeta$, we must treat the conformal and nonconformal cases separately.

\textbf{Assumption:} We assume that $P(0)\neq 0$, which corresponds to a nonconformal harmonic map. This assumption will persist until Section~\ref{sec:conformal}, where we deal with the conformal case.

Having made this assumption $\zeta$ is not a factor of $P$ and so cannot be a factor of $FF^i$. Therefore $FF^i$ divides $c^i$. Applying this to~\eqref{eqn:Q},
\[
FF^1G\tilde{b}^1 c^2 - FF^2G\tilde{b}^2 c^1 = Q FF^1F^2\tilde{P}.
\]
By definition, neither $F$ nor $G$ divide $\tilde{P}$, demonstrating that $FG$ divides $Q$. This provides a bound on the number of coincident roots that are allowed if a deformation is to exist; $Q$ is quadratic so $FG = \gcd(b^1,b^2)$ must be degree two or less. Moreover, because all of $P$, $b^1$, $b^2$ are real, and $P$ has no roots on the unit circle, any common roots of the three polynomials must come in conjugate inverse pairs and so the degree of $F$ will always be even.

Assuming that a deformation does exist, $FF^i$ must divide $c^i$ and $FG$ must divide $Q$. 
If the polynomials $c^i$ or $Q$ did not have the factors indicated by~\eqref{eqn:def c Q tilde}, then there would be factors on the left hand sides of~\eqref{eqn:EMPDi} and~\eqref{eqn:Q} that did not appear on the right hand sides, and this contradiction would preclude the possibility of a solution.
Hence equations~\eqref{eqn:EMPDi reduced} and~\eqref{eqn:Q reduced} are necessary conditions to be able to solve~\eqref{eqn:EMPDi} and~\eqref{eqn:Q} respectively. 
Let us define
\[
c^i = FF^i \tilde{c}^i, \text{ and } Q = FG\tilde{Q}.
\labelthis{eqn:def c Q tilde}
\]
We can then remove the common factor $FF^i$ from~\eqref{eqn:EMPDi} to arrive at the reduced equations
\[
\dot{P} G \tilde{b}^i - 2 F^j \tilde{P} \dot{b}^i = 2 F^j \tilde{P} (\hat{c}^i - \zeta\hat{c}^{i\prime}) + \zeta(\zeta^2-1)P' \tilde{c}^i,
\labelthis{eqn:EMPDi reduced}
\]
for $i=1,2$ and $j\neq i$. In the same manner, the $Q$ equation~\eqref{eqn:Q} reduces to
\begin{align*}
FF^1G\tilde{b}^1 FF^2\tilde{c}^2 - FF^2G\tilde{b}^2 FF^1\tilde{c}^1 &= FG\tilde{Q} FF^1F^2\tilde{P} \\
\tilde{b}^1 \tilde{c}^2 - \tilde{b}^2 \tilde{c}^1 &= \tilde{Q} \tilde{P}.
\labelthis{eqn:Q reduced}
\end{align*}

To recap, any tangent vector $(\dot P, \dot b^1, \dot b^2)$ to $\mathcal{M}_g$ must satisfy~\eqref{eqn:residueTangent}. Each tangent vector gives rise to a pair of polynomials $(c^1,c^2)\in\mathcal{P}^{g+1}_\R \times \mathcal{P}^{g+1}_\R$ through~\eqref{eqn:EMPDi}. The polynomials $(c^1,c^2)$ must themselves factorise as per~\eqref{eqn:def c Q tilde} and satisfy~\eqref{eqn:Q reduced} for some real polynomial $\tilde{Q}$ of degree at most $2$. These are necessary and sufficient conditions for a triple of polynomials to be a tangent vector to $\mathcal{M}_g$. Under various additional assumptions, Lemmata~\ref{lem:tangent generic},~\ref{lem:tangent G}, and~\ref{lem:tangent conformal} find all solutions to these equations.

%%%%%%%%%%%%%%%%%%%%%%%%%%%%%%%%%%%%%%%%%%%%%%%%%%%%%%%%%%%%%%%%%%%%
%%%%%%%%%%%%%%%%%%%%%%%%%%%%%%%%%%%%%%%%%%%%%%%%%%%%%%%%%%%%%%%%%%%%
%%%%%%%%%%%%%%%%%%%%%%%%%%%%%%%%%%%%%%%%%%%%%%%%%%%%%%%%%%%%%%%%%%%%
%%%%%%%%%%%%%%%%%%%%%%%%%%%%%%%%%%%%%%%%%%%%%%%%%%%%%%%%%%%%%%%%%%%%
\section{B\'ezout Identity}\label{sec:bezout identity}
Equations such as~\eqref{eqn:EMPDi} and~\eqref{eqn:Q reduced} are of the form
\[
AX - BY = C,
\labelthis{eqn:ABC}
\]
to which B\'ezout's identity for polynomials applies. This section will develop a variant of B\'ezout's identity adapted to the particulars of our situation. The basic version of B\'ezout's identity for polynomials asserts that if $\gcd(A,B) = 1$, then there is a unique solution $(\mathbf{X},\mathbf{Y})$ of minimal degree, with $\deg \mathbf{X} < \deg B$ and $\deg \mathbf{Y} < \deg A$. Let us fix $A, B, C$ to be polynomials of degree $a,b, c$ respectively and $d$ to be the degree of $D := \gcd(A,B)$.

\textbf{Assumption:} Throughout this section, we assume that $\gcd(A,B)$ divides $C$. This is a necessary condition to the solution of~\eqref{eqn:ABC}.

\begin{lem}\label{lem:ABC soln}
There is a unique solution $(\mathbf{X},\mathbf{Y})$ to~\eqref{eqn:ABC} such that the degree of $\mathbf{X}$ is at most $b-d-1$. Moreover, if $c < a+b-d$ then the degree of $\mathbf{Y}$ is at most $a-d-1$.

\begin{proof}
% One may prove the existence of these minimal solutions using the Euclidean algorithm~\cite{Mora2003}, but 
Because we will have need of the specific formula~\eqref{eqn:linear system} we construct the minimal solution by giving a linear system of equations that the coefficients of $\mathbf{X}$ must satisfy and showing that there is a unique solution to this linear system.

Suppose then that $\mathbf{X}$ is degree $b-d-1$. By dividing~\eqref{eqn:ABC} by the common factor $D = \gcd(A,B)$, consider the equation
$(A/D) X - (B/D) Y = C/D$.
If $\beta$ is a root of $B/D$ then evaluation yields,
\[
X_0 + X_1 \beta + X_2 \beta^2 + \cdots + X_{b-d-1} \beta^{b-d-1} = \frac{(C/D)(\beta)}{(A/D)(\beta)} = (C/A)(\beta).
\]
If $\beta$ is a root of $B$ of multiplicity $r$ then we may differentiate repeatedly to obtain $r$ linearly independent equations. If we label the distinct roots of $B/D$ as $\beta_i$, their multiplicities as $r_i$, and for brevity let $n=b-d-1$ then the full system of equations is
% to obtain another relation,
% \begin{align*}
% X &= \frac{C}{A} - B\frac{Y}{A} \\
% X' &= \bra{\frac{C}{A}}' - B'\frac{Y}{A} - B\bra{\frac{Y}{A}}'\\
% \mathbf{X}'(\beta) &= X_1 + 2 X_2 \beta + \cdots + (b-d-1) X_{b-d-1} \beta^{b-d-2} = \bra{\frac{C}{A}}'(\beta),
% \end{align*}
% Differentiating 
\begin{longeqn}
\begin{bmatrix}
  1 & \beta_1 & (\beta_1)^2 & \ldots & (\beta_1)^{r_1-1} & \ldots & (\beta_1)^{n} \\
  0 & 1 & 2\beta_1 & \ldots & (r_1-1)(\beta_1)^{r_1-2} & \ldots & n(\beta_1)^{n-1} \\
  \vdots & \vdots & \vdots & & \vdots & & \vdots \\
  0 & 0 & 0 & \ldots & 1 & \ldots & \frac{n!}{(n+1 - r_1)!}(\beta_1)^{n+1-r_1} \\
  1 & \beta_2 & (\beta_2)^2 & \ldots & (\beta_2)^{r_2-1} & \ldots & (\beta_2)^{n} \\
%   \vdots & \vdots & \vdots & & \vdots & & \vdots \\
%   0 & 0 & 0 & \ldots & 1 & \ldots & \frac{n!}{(n+1 - r_2)!}(\beta_2)^{n+1-r_2} \\
%   1 & \beta_3 & (\beta_3)^2 & \ldots & (\beta_3)^{r_2-1} & \ldots & (\beta_3)^{n} \\
  \vdots & \vdots & \vdots & & \vdots & & \vdots \\
  \vdots & \vdots & \vdots & & \vdots & & \vdots \\
  0 & 0 & 0 & \ldots & 1 & \ldots & \frac{n!}{(n+1 - r_k)!}(\beta_k)^{n+1-r_k} \\
\end{bmatrix}
\begin{bmatrix}
  X_0 \\ X_1 \\ \vdots \\ X_n
\end{bmatrix}
=
\begin{bmatrix}
  (C/A)(\beta_1) \\ (C/A)'(\beta_1) \\ \vdots \\ (C/A)^{(r_1-1)}(\beta_1) \\
  (C/A)(\beta_2) \\ 
%   \vdots \\ (C/A)^{(r_2-1)}(\beta_2) \\
%   (C/A)(\beta_3) 
  \\ \vdots \\ \vdots \\ (C/A)^{(r_k-1)}(\beta_k) \\
\end{bmatrix},
\end{longeqn}
\[
\text{     }\labelthis{eqn:linear system}
\]

The $(n+1)\times (n+1)$ coefficient matrix on the left is called the confluent Vandermonde matrix at the roots of $B/D$, and we shall denote it $V(B/D)$. We shall denote the vector on the right by $h(B/D,C/A)$. A confluent Vandermonde matrix is always nonsingular~\cite{Kalman1984}, therefore there is a unique solution to this system.

Having found $\mathbf{X}$, we may find the corresponding $\mathbf{Y}$ from $B\mathbf{Y} = A\mathbf{X} - C$. To address the final part of the lemma, if $c < a+b-d$, then the degree of the right-hand side above is degree at most $a+b-d-1$, whence the degree of $\mathbf{Y}$ is degree at most $a-d-1$.
\end{proof}
\end{lem}

\begin{lem}\label{lem:ABC real soln}
Suppose that $(X,Y)$ is a solution to equation~\eqref{eqn:ABC}.
Further suppose that $A,B$ and $C$ are members of $\mathcal{P}^a_\R,\mathcal{P}^b_\R$ and $\mathcal{P}^c_\R$ respectively, and that $\deg X \leq c-a$.
Then there exists a solution to the equation belonging to $\mathcal{P}^{c-a}_\R \times \mathcal{P}^{c-b}_\R$.
Moreover, if $c < a+b-d$ then the minimal solution $(\mathbf{X},\mathbf{Y})$ itself belongs to $\mathcal{P}^{c-a}_\R \times \mathcal{P}^{c-b}_\R$.
\begin{proof}
% Recall that the polynomials of $\mathcal{P}^a_\R$ are exactly those that are fixed under the real involution $\rho_a$. 
% Explicitly, for $A \in \mathcal{P}_\R^a$
% \[
% A(\zeta) = \rho_a^* A := \bar{\zeta}^a \bar{A}\left(\zeta^{-1}\right),
% \]
% where $\bar{A}(\zeta)$ is the polynomial obtained from $A(\zeta)$ by conjugating its coefficients. 
% Every polynomial $A$ of degree not exceeding $a$ may be written as the sum of its real and imaginary parts
% \[
% A = \frac{1}{2}\left( A + \rho_a^* A \right) + \frac{1}{2}\left( A - \rho_a^* A \right),
% \]
% and it is an easy check that a polynomial belongs to $\mathcal{P}_\R^a$ iff the imaginary part vanishes.

Suppose that the equation $AX-BY=C$ has a solution $(X,Y)$. By applying the involution of $\mathcal{P}_\R^c$ to~\eqref{eqn:ABC} we see that
\begin{align*}
% \bar{\zeta}^c\bar{A}(\zeta^{-1})\bar{X}(\zeta^{-1})
% -\bar{\zeta}^c\bar{B}(\zeta^{-1})\bar{Y}(\zeta^{-1})
% &= \bar{\zeta}^c\bar{C}(\zeta^{-1}) \\
%%%%%%%%%%%%%%%%%%%%%%
A(\zeta)\bar{\zeta}^{c-a}\bar{X}(\zeta^{-1})
-B(\zeta)\bar{\zeta}^{c-b}\bar{Y}(\zeta^{-1})
&= C(\zeta).
\labelthis{eqn:involution_ABC}
\end{align*}
By assumption, the degree of $X$ is at most $c-a$ so $\bar{\zeta}^{c-a}\bar{X}(\zeta^{-1})$ is polynomial. From $AX - BY = C$ it follows that $Y$ is degree at most $c-b$ and hence $\bar{\zeta}^{c-b}\bar{Y}(\zeta^{-1})$ is a polynomial too. By averaging equation~\eqref{eqn:involution_ABC} with the original equation~\eqref{eqn:ABC}, we arrive at the following,
\[
A\frac{1}{2}\left(X + \rho_{c-a}^*X \right)
-B\frac{1}{2}\left(Y + \rho_{c-b}^*Y \right)
= C,
\]
yielding a real solution as claimed.

In the case that we also have $c<a+b-d$, by subtracting~\eqref{eqn:involution_ABC} from the original equation~\eqref{eqn:ABC}, we arrive at the following equation for the `imaginary parts'
\begin{align*}
% A\left(\mathbf{X} - \rho_{c-a}^*\mathbf{X} \right)
% -B\left(\mathbf{Y} - \rho_{c-b}^*\mathbf{Y} \right)
% &= 0 \\
%%%%%%%%%%%%%%%%%%%%%%%%%%%%
(A/D)\left(\mathbf{X} - \rho_{c-a}^*\mathbf{X} \right)
-(B/D)\left(\mathbf{Y} - \rho_{c-b}^*\mathbf{Y} \right)
&= 0.
\end{align*}
We see therefore that $B/D$ must divide $\mathbf{X} - \rho_{c-a}^*\mathbf{X}$, but the former is degree $b-d$ whereas the latter is at most degree $d-b-1$. It follows therefore that $\mathbf{X}=\rho_{c-a}^*\mathbf{X}$ and so $\mathbf{X} \in \mathcal{P}_\R^{c-a}$. Similarly we have that $\mathbf{Y} \in \mathcal{P}_\R^{c-b}$.
\end{proof}
\end{lem}

\begin{lem}\label{lem:ABC soln space}
Suppose that the polynomials $A,B,C$ are members of $\mathcal{P}^a_\R, \mathcal{P}^b_\R, \mathcal{P}^c_\R$ respectively. If $c \geq a+b-d$, then the space of solutions $(X,Y) \in \mathcal{P}^{c-a}_\R \times \mathcal{P}^{c-b}_\R$ to equation~\eqref{eqn:ABC} is
\[
\left\{(\mathbf{X} + U(B/D), \mathbf{Y} + U(A/D)) \mid U \in \mathcal{P}_R^{c-a-b+d} \right\},
\]
where $(\mathbf{X},\mathbf{Y})$ is the minimal solution to~\eqref{eqn:ABC} given in Lemma~\ref{lem:ABC soln}.

\begin{proof}
First suppose that $(Z,W)$ is a solution to the related equation $AZ-BW = 0$. After removing the common factor $D$, it must be that $B/D$ divides $Z$. For some polynomial $U$, let $Z = (B/D)U$. Then
\begin{align*}
(A/D)(B/D)U - (B/D) W = 0,
\end{align*}
from which we conclude that $W = (A/D)U$. If $Z$ and $W$ are real, so too is the quotient $U$.

Return now to the equation $AX-BY=C$ and suppose then that $(\mathbf{X},\mathbf{Y})$ the minimal solution given by Lemma~\ref{lem:ABC soln}.
Because $c-a \geq b-d > \deg \mathbf{X}$, by Lemma~\ref{lem:ABC real soln} this is a real solution in $\mathcal{P}^{c-a}_R \times \mathcal{P}^{c-b}_R$.

If $(X,Y)$ is any other real solution, then $A(X-\mathbf{X}) - B(Y-\mathbf{Y}) = 0$. Therefore $X - \mathbf{X} = (B/D)U$ and $Y - \mathbf{Y} = (A/D)U$ for some real polynomial $U$ in $\mathcal{P}_\R^{c-a-b+d}$. Conversely, given any solution $(X,Y)$, it is clear that $(X+(B/D)U,Y+(A/D)U)$ is again a real solution for any real polynomial $U$.
Thus polynomials of this form are exactly the desired solutions to $AX-BY=C$.
\end{proof}
\end{lem}

%%%%%%%%%%%%%%%%%%%%%%%%%%%%%%%%%%%%%%%%%%%%%%%%%%%%%%%%%%%%%%%%%%%%
%%%%%%%%%%%%%%%%%%%%%%%%%%%%%%%%%%%%%%%%%%%%%%%%%%%%%%%%%%%%%%%%%%%%
%%%%%%%%%%%%%%%%%%%%%%%%%%%%%%%%%%%%%%%%%%%%%%%%%%%%%%%%%%%%%%%%%%%%
%%%%%%%%%%%%%%%%%%%%%%%%%%%%%%%%%%%%%%%%%%%%%%%%%%%%%%%%%%%%%%%%%%%%
\section{The Tangent Space to \texorpdfstring{$\mathcal{M}_g$}{Mg}}\label{sec:tangent space}

In Section~\ref{sec:deformations}, we elucidated several necessary properties of an infinitesimal deformation of spectral triples. Now we turn our attention to the converse; in Lemmata~\ref{lem:tangent generic} and~\ref{lem:tangent G} we give conditions under which it is possible to solve~\eqref{eqn:EMPDi} and~\eqref{eqn:Q reduced} and thereby find an infinitesimal deformation. Firstly, we examine whether it is possible to construct polynomials $c^i$ that factor as per~\eqref{eqn:def c Q tilde} and solve~\eqref{eqn:Q reduced} for a given $\tilde{Q}$, and whether this construction is unique. Secondly, we shall insert the polynomials $c^i$ into the right hand side of~\eqref{eqn:EMPDi} and solve it to recover $\dot{b}^i$ and $\dot{P}$.

For each equation, the main obstacle to the existence of a solution is common factors among the polynomials $(P,b^1,b^2)$. If there are too many common factors (for example if $\gcd(b^1,b^2)$ has degree greater than two), then it will not be possible to deform the spectral data. Even when it is possible to deform, the form of the solution of~\eqref{eqn:Q reduced} is dependent on those common factors. Thus we will need to divide our approach into several cases according to the common roots of $b^1$ and $b^2$, labelled in the following table. Refer to Figure~\ref{fig:common factors} for the definitions of $F$ and $G$.
\begin{center}
\begin{tabular}{|c|c|c|c|}
\hline
Case & $P_0$ & $\deg F$ & $\deg G$ \\ \hline\hline
(a) & \multirow{4}{*}{$P_0 \neq 0$} & 0 & 0 \\ \cline{1-1}\cline{3-4}
(b) && 0 & 1, 2 \\ \cline{1-1}\cline{3-4}
(c) && 2 & 0 \\ \cline{1-1}\cline{3-4}
(d) && \multicolumn{2}{|c|}{$FG \in \mathcal{P}^k_\R,\, k>2$} \\ \hline
\end{tabular}
\end{center}
Recall our standing assumptions that the spectral curve is nonsingular and that it is not branched over $\zeta=0$, ie that the spectral triple is from a nonconformal harmonic map. This is equivalent to $P_0 \neq 0$. These four cases are disjoint and exhaustive for nonconformal spectral triples. 

\begin{defn}\label{def:subsets U}
We denote the points of case (a) in $\mathcal{U}$ as $\mathcal{U}_{a}$, and likewise for the other cases. Let $\mathcal{U}_{ab} = \mathcal{U}_{a}\cup \mathcal{U}_{b}$. We may equivalently characterise $\mathcal{U}_{ab}$ as the set
\[
\mathcal{U}_{ab} = \Set{ (P,b^1,b^2) \in \mathcal{U} } { \gcd(P,b^1,b^2) = 1 \text{ and } \gcd(b^1,b^2) \in \mathcal{P}^l_\R \text{ for } l \leq 2 },
\]
which is open and dense in $\mathcal{U}$.
\end{defn}

% They may also be described purely in terms of the greatest common factor of $b^1$ and $b^2$.
% \begin{center}
% \begin{tabular}{|c|l|}
% \hline
% Case & Description\\ \hline\hline
% a & $\gcd(b^1,b^2) = 1$\\ \hline
% b & $\gcd(b^1,b^2)$ is linear or quadratic and does not divide $P$ \\ \hline
% c & $\gcd(b^1,b^2)$ is quadratic and divides $P$ \\ \hline
% d & The degree of $\gcd(b^1,b^2)$ is greater than two \\ \hline
% \end{tabular}
% \end{center}

The remainder of the section proves that $\mathcal{M}_g \cap \mathcal{U}_{ab}$ is a manifold. We have seen that the spectral triples in case (d) do not admit any deformations. In general, points of $\mathcal{M}_g$ where case (c) holds are singularities of a deformation~\cite{Schmidt2016} and are not considered further in this paper.

% We may characterise $\mathcal{U}_{ab}$ as the set
% \[
% \mathcal{U}_{ab} = \Set{ (P,b^1,b^2) \in \mathcal{U} } { \gcd(P,b^1,b^2) = 1 \text{ and } \gcd(b^1,b^2) \in \mathcal{P}^l_\R \text{ for } l \leq 2 }.
% \]
% It is open because the complement $\mathcal{U}\setminus \mathcal{U}_{ab}$ may be written as the union of three closed sets:
% \begin{align*}
% &\Set{(P,b^1,b^2) \in \mathcal{U}}{ \gcd(b^1,b^2) \in\mathcal{P}^k_\R,\,k > 2 } \\
% &\qquad\cup
% \Set{(P,b^1,b^2) \in \mathcal{U}}{ \gcd(P,b^1,b^2) \in\mathcal{P}^l_\R,\,l \geq 2 } \\
% &\qquad\cup
% \Set{(P,b^1,b^2) \in \mathcal{U}}{ P_0 = 0 }.
% \end{align*}

Before we proceed to solving equations, there is a short digression we must make. In case (a) given an arbitrary quadratic polynomial $Q$ it will not always be possible to solve~\eqref{eqn:Q reduced} for polynomials $\tilde{c}^i(\zeta)$ corresponding to an infinitesimal deformation $(\dot{P},\dot{b}^1,\dot{b}^2)$. To see why, consider the linear system of equations in the coefficients of $\tilde{c}^2$ that arises from evaluating~\eqref{eqn:Q reduced} at the roots of $\tilde{b}^2$,
\[
V\left(\tilde{b}^2\right)
\begin{bmatrix}
\tilde{c}^2_0 \\ \vdots \\ \tilde{c}^2_n
\end{bmatrix}
=
h\bra{\tilde{b}^2, \frac{Q\tilde{P}}{\tilde{b}^1} },
\labelthis{eqn:Q reduced linear}
\]
where $V$ is the confluent Vandermonde matrix defined by~\eqref{eqn:linear system}. 
We know that the degree of the unique minimal solution could be as high as
\[
n := \deg \tilde{b}^2 - 1 = g + 2 - d_2.
\]
But 
% $c^2$ is degree $g+1$, and $c^2 = F^2\tilde{c}^2$. Therefore 
the degree of $\tilde{c}^2$ is $g + 1 - d_2 = n-1$. 
% If every solution of~\eqref{eqn:Q reduced} is degree $n$ or more, then there can be no solutions that correspond to infinitesimal deformations. 
Thus we must introduce a restriction on our choice of $Q$ so that a solution to~\eqref{eqn:Q reduced} of the correct degree exists. We will express this restriction in terms of the vanishing of a function $R$.

\begin{defn}\label{def:def R}
Recall the confluent Vandermonde matrix $V(B)$ and vector $h(B,C/A)$ defined in~\eqref{eqn:linear system}. We define the function $R$ to be
\begin{align*}
R : \mathcal{U}_{a} \times \mathcal{P}^2_\R &\to \C \\
((P,b^1,b^2), Q) &\mapsto \text{ the last entry of } \left[V( \tilde{b}^2 )\right]^{-1}
h\bra{\tilde{b}^2, \frac{Q\tilde{P}}{\tilde{b}^1} }.
\end{align*}
When the point of $\mathcal{U}_{a}$ is understood, we shall abbreviate this to $R(Q)$.
\end{defn}

This function $R$ is simply the function that gives the value of the degree $n$ coefficient of $\tilde{c}^2$; the condition that $R((P,b^1,b^2),Q) = 0$ is equivalent to the condition that there is a solution $\tilde{c}^2$ to~\eqref{eqn:Q reduced} of degree $n-1$ or less.
Likewise, evaluating~\eqref{eqn:Q reduced} at the roots of $\tilde{b}^1$ leads to a solution $\tilde{c}^1$ of degree $g+2 - d_1$. From the highest order term of~\eqref{eqn:Q reduced}, if $\tilde{c}^2$ is degree $n-1$ or lower, then $\tilde{c}^1$ will be degree $g + 1 - d_1$ or lower without any further restrictions on $Q$.

In case (b), the degree of $\tilde{b}^2$ is
\[
g+3 - (d_2 + d_F + d_G) \leq g + 2 - d_2,
\]
as $G$ will be nontrivial. This means that the minimal solution to~\eqref{eqn:Q reduced}, which has degree strictly less than $\deg \tilde{b}^2$, will have the correct degree without needing to impose any extra conditions on $Q$.

It is important to note that $R$ is a linear function in the coefficients of $Q$. $R$ satisfies the following reality type condition and thus at any point of $\mathcal{U}_{a}$ there is a real 2-plane of polynomials $Q\in\mathcal{P}^2_\R$ that satisfy $R(Q) = 0$.

\begin{lem}
At every point $(P,b^1,b^2)$ of $\mathcal{U}_{a}$, $R((P,b^1,b^2), Q)$ satisfies the relation
\[
\bar{R} = (-1)^{n}\bra{ \prod_{i=1}^{n+1}  \beta_i }  R,
\labelthis{eqn:R reality}
\]
where $\beta_i$ are the $n+1 = g+3 - d_2$ roots of $\tilde{b}^2$, counted with multiplicity.

\begin{proof}
We shall demonstrate this property first at points where the roots of $\tilde{b}^2$ are distinct. Let $\tilde{b}^2$ have $n+1 = g+3 - d_2$ distinct roots $\beta_i$. In this case, the explicit form of the solution to the linear system of equations~\eqref{eqn:Q reduced linear} is elegant. Consider the Lagrange polynomials at the roots of $\tilde{b}^2$,
\[
L_i (\zeta) := \prod_{j \neq i} \frac{\zeta-\beta_j}{\beta_i - \beta_j}.
\]
Each of these polynomials is degree $n$ and has the property that $L_i (\beta_j) = \delta_{ij}$. The unique polynomial of degree at most $n$ solving the linear system is
\[
\tilde{c}^2(\zeta) = \sum_{i = 1}^{n+1} \bra{\frac{Q\tilde{P}}{\tilde{b}^1}}(\beta_{i}) L_i (\zeta),
\]
and in particular the highest coefficient is $R$,
\[
R = \sum_{i = 1}^{n+1} \bra{\frac{Q\tilde{P}}{\tilde{b}^1}}(\beta_i) \prod_{j \neq i} \bra{ \beta_i - \beta_j }^{-1}.
\]
Because $\tilde{b}^2$ is a real polynomial, its set of roots is invariant under $\zeta \mapsto \cji{\zeta}$. This creates an involution on the set of roots. Let $\tau$ be the involution on the integers $\{1,2,\ldots,n+1\}$ such that $\beta_{\tau(i)} = \cji{\beta}_i$. We compute that
% \[
% \bar{ \bra{\frac{Q\tilde{P}}{\tilde{b}^1}}(\beta_{i}) }
% = \bar{\frac{Q(\beta_i)\tilde{P}(\beta_i)}{\tilde{b}^1(\beta_i)}}
% = \bar{\beta}_i^{g+1 - d_2} \frac{Q(\cji{\beta}_i)\tilde{P}(\cji{\beta}_i)}{\tilde{b}^1(\cji{\beta}_i)}
% = \bar{\beta}_i^{n-1} \bra{\frac{Q\tilde{P}}{\tilde{b}^1}}(\beta_{\tau(i)}),
% \]
% and
% \begin{align*}
% \prod_{j \neq i} \bra{ \bar{\beta}_i - \bar{\beta}_j }^{-1}
% &= \prod_{j \neq i} \cji{\beta}_i \cji{\beta}_j \bra{ \cji{\beta}_j - \cji{\beta}_i }^{-1} \\
% &= \bra { \cji{\beta}_i }^{n} \bra{ \prod_{j \neq i}  \beta_{\tau(j)}} (-1)^{n} \prod_{j \neq i}\bra{ \beta_{\tau(i)} - \beta_{\tau(j)} }^{-1},
% \end{align*}
% so that the conjugate of $R$ is
\begin{align*}
\bar{R}
% &= \sum_{i = 1}^{n+1} \bra{\frac{Q\tilde{P}}{\tilde{b}^1}}(\beta_{i}) \prod_{j \neq i} \bra{ \bar{\beta}_i - \bar{\beta}_j }^{-1} \\
%%%%%%%%%%%%%%%%%%%%%%%%%%%%%%%%%%%%
&= \sum_{i = 1}^{n+1} 
\left[\rule{0cm}{0.7cm}  \bra{\frac{Q\tilde{P}}{\tilde{b}^1}}(\beta_{\tau(i)})\,{\bar{\beta}}_i^{n-1} \right]
\left[\rule{0cm}{0.7cm}  \prod_{j \neq i} \cji{\beta}_i \cji{\beta}_j \bra{ \beta_{\tau(j)} - \beta_{\tau(i)} }^{-1}\right] \\
%%%%%%%%%%%%%%%%%%%%%%%%%%%%%%%%%%%%
&= (-1)^{n}\bra{ \prod_{j=1}^{n+1}  \beta_j }  R.
\end{align*}

To complete the proof, we must extend the argument to points $p = (P,b^1,b^2)$ where the roots of $\tilde{b}^2$ are not distinct.
% we may construct a sequence of points $p_k$ of $\mathcal{U}_{a}$ that converge to the point $p$, but such that $b^2(p_k)(\zeta)$ has distinct roots for every $k$. We will show that $\lim_{k\to\infty} R(p_k,Q) = R(p,Q)$, and therefore the relation holds in the limit.
% \begin{align*}
% \bar{R(p,Q)}
% = \lim_{k\to\infty} \bar{R(p_k,Q)}
% &= \lim_{k\to\infty} (-1)^{n}\bra{ \prod_{i=1}^{n+1}  \beta_i(p_k) }  R(p_k,Q) \\
% &= (-1)^{n}\bra{ \prod_{i=1}^{n+1}  \beta_i } R(p,Q).
% \end{align*}
Suppose that we are at a point $p = (P,b^1,b^2)$ of $\mathcal{U}_{a}$ where $\tilde{b}^2$ has a double root $\beta$. Considering the subvariety of $\mathcal{U}_{a}$ where $F^2 = \gcd(P,b^2)$ is fixed, we may find a sequence of points $p_k$ in this subvariety converging to $p$ with the property that the roots of each polynomial $\tilde{b}^2(p_k)$ are distinct. Let us label the two simple roots of $\tilde{b}^2(p_k)$ that coalesce at $p$ to form the double root $\beta$ as $\beta_1(k),\beta_2(k)$.
% In other words, these are the two roots of $\tilde{b}^2(p_k)$ such that $\lim_{k\to\infty} \beta_1(p_k) = \lim_{k\to\infty} \beta_2(p_k) = \beta$. The corresponding rows of the Vandermonde matrix $V(\tilde{b}^2(p_k))$ in the system~\eqref{eqn:Q reduced linear} are
% \[
% \begin{bmatrix}
% 1 & \beta_1 & (\beta_1)^2 & \cdots & (\beta_1)^{n} \\
% 1 & \beta_2 & (\beta_2)^2 & \cdots & (\beta_2)^{n} \\
% &&&\vdots&
% \end{bmatrix}
% \begin{bmatrix}
% \tilde{c}^2_0 \\~\\ \tilde{c}^2_1 \\~\\ \vdots \\~\\ \tilde{c}^2_n
% \end{bmatrix}
% =
% \begin{bmatrix}
% \bra{\frac{Q\tilde{P}}{\tilde{b}^1}}(\beta_1) \\
% \bra{\frac{Q\tilde{P}}{\tilde{b}^1}}(\beta_2) \\
% \vdots
% \end{bmatrix}.
% \]

Consider the corresponding rows of the Vandermonde matrix $V(\tilde{b}^2(p_k))$. Performing elementary row operations does not change the solution to this system, and so we may subtract the one row from the other and scale it by $(\beta_2-\beta_1)^{-1}$. This gives
\begin{longeqn}
\begin{bmatrix}
    1 & \beta_1 & (\beta_1)^2  & \ldots & (\beta_1)^{n} \\
    0 & 1 & \beta_2+\beta_1  & \ldots & \sum_{j=0}^{n-1}(\beta_1)^j(\beta_2)^{n-1-j} \\
    &&&\vdots&
\end{bmatrix}
\begin{bmatrix}
    \tilde{c}^2_0 \\ \vdots \\ \tilde{c}^2_n
\end{bmatrix}
=
\begin{bmatrix}
    \bra{\frac{Q\tilde{P}}{\tilde{b}^1}}(\beta_1) \\
    (\beta_2 - \beta_1)^{-1}\left[ \bra{\frac{Q\tilde{P}}{\tilde{b}^1}}(\beta_2) - \bra{\frac{Q\tilde{P}}{\tilde{b}^1}}(\beta_1) \right] \\
    \vdots
    \end{bmatrix}.
\end{longeqn}
The limit of the above as $k \to \infty$ is precisely the confluent Vandermonde matrix at $\tilde{b}^2(p)$. The calculation for higher order roots is similar. If $\tilde{b}^2$ has more than one higher order root, then we may perform this operation concurrently for each of them.

Combining these row operations with the fact that inversion of a matrix is continuous and that the roots of a polynomial are continuous functions of coefficients~\cite[Theorem V.4A]{Whitney1972}, this shows that limit of the solutions $\tilde{c}^2(p_k)$ is just the solution $\tilde{c}^2(p)$. In particular, the last component of $\tilde{c}^2(p)$ is $R(p,Q)$, and so $\lim_{k\to\infty} R(p_k,Q) = R(p,Q)$. Hence we have established~\eqref{eqn:R reality} at all points of $\mathcal{U}_{a}$.
\end{proof}
\end{lem}

At this point we are ready to solve equations~\eqref{eqn:EMPDi} and~\eqref{eqn:Q reduced} in both of the cases (a) and (b). At the end of that process, we will have constructed a tangent vector $(\dot{P}, \dot{b}^1, \dot{b}^2)$ to the space of spectral triples $\mathcal{M}_g$.

\begin{lem}[Case (a)]\label{lem:tangent generic}
Take a spectral triple $(P,b^1,b^2)\in\mathcal{M}_g$ associated with a nonconformal harmonic map, with a nonsingular spectral curve given by $\eta^2 = P$ of genus $g$. Suppose that $\gcd(b^1,b^2) = 1$. Then for every polynomial $Q \in \mathcal{P}^2_\R$ with $R(Q) = 0$, there exist unique real polynomials $c^i\in\mathcal{P}^{g+1}_\R$ that factor as per~\eqref{eqn:def c Q tilde} and satisfy~\eqref{eqn:Q reduced}. Further, for each such pair $(c^1,c^2)$, there is a unique vector $(\dot P, \dot b^1, \dot b^2) \in \mathcal{P}_\R^{2g+2}\times\mathcal{P}_\R^{g+3}\times\mathcal{P}_\R^{g+3}$ that satisfies~\eqref{eqn:EMPDi} and~\eqref{eqn:residueTangent}. It is therefore a tangent vector to the space of spectral triples $\mathcal{M}_g$.

\begin{proof}
Let us first begin with~\eqref{eqn:Q reduced}. By Lemma~\ref{lem:ABC soln} there is a unique solution $(\mathbf{\tilde{c}}^1, \mathbf{\tilde{c}}^2)$ to this equation of degree at most $(g+2-d_1, g+2-d_2)$, where $d_i = \deg F^i$. We note that the leading coefficient $\mathbf{\tilde{c}}^2_{g+2-d_2}$ is $R(Q)$ by definition, which by assumption is zero.
Examining the highest order of~\eqref{eqn:Q reduced}, if $\mathbf{\tilde{c}}^2_{g+2-d_2}$ vanishes, so too must $\mathbf{\tilde{c}}^1_{g+2-d_1}$. Multiplying~\eqref{eqn:Q reduced} through by $F^1F^2$, we arrive at unique polynomials $c^i = F^i \mathbf{\tilde{c}}^i$ with the factors required by~\eqref{eqn:def c Q tilde}. Both of these polynomials are degree at most $g+1$. By Lemma~\ref{lem:ABC real soln} these are real polynomials. We similarly define $\hat{c}^i = (\zeta^2 -1)c^i$.

Next we must solve~\eqref{eqn:EMPDi}. To see that Lemma~\ref{lem:ABC real soln} will apply, for any $f\in\mathcal{P}_\R^{k}$ we compute that
\[
\rho^*_k(\zeta f') = \zeta^{k-1}\bar{f}'(\zeta^{-1}) = k f(\zeta) - \zeta f'(\zeta).
\]
From this it follows that the right hand side of~\eqref{eqn:EMPDi} is real for all real polynomials $P$ and imaginary polynomials $\hat{c}^i$, not just when these polynomials arise from a deformation.
% Equation~\eqref{eqn:EMPDi reduced} is a necessary condition for~\eqref{eqn:EMPDi}. Equation~\eqref{eqn:EMPDi reduced} in this case it reads
% \[
% \dot{P} \tilde{b}^i - 2 F^j \tilde{P} \dot{b}^i = 2 F^j \tilde{P} (\hat{c}^i - \zeta\hat{c}^{i\prime}) + \zeta(\zeta^2-1)P' \tilde{c}^i.
% \labelthis{eqn:EMPDi nc FG=1}
% \]
It can be solved using a Vandermonde matrix, if $\tilde{b}^i$ is non-vanishing at the roots of $F^j \tilde{P}$ and vice versa. But by definition they are coprime.

The two equations for $i=1,2$ may give different solutions for $\dot P$, and indeed in general they do. However, we will show that there is a common solution to both. Let a solution to each equation~\eqref{eqn:EMPDi} be $(\mathbf{\dot{P}}^1, \mathbf{\dot{b}}^1)$ and $(\mathbf{\dot{P}}^2, \mathbf{\dot{b}}^2)$. From Lemma~\ref{lem:ABC soln space} the sets of solutions of degree $(2g+2,g+3)$ are
\[
\Set { \bra{ \mathbf{\dot{P}}^1 + 2rF^2\tilde{P}, \mathbf{\dot{b}}^1 + r\tilde{b}^1 } }
{ r \in\mathcal{P}^{d_1}_\R },
\]
and
\[
\Set { \bra{ \mathbf{\dot{P}}^2 + 2sF^1\tilde{P}, \mathbf{\dot{b}}^2 + s\tilde{b}^2 } }
{ s \in\mathcal{P}^{d_2}_\R},
\]
respectively. First note that every element of both of these sets take the same value at any root $\alpha$ of $\tilde{P}$. This follows from
\[
\dot P^1(\alpha)
= \alpha (\alpha^2 - 1) P'(\alpha) \frac{\tilde{c}^1(\alpha)}{\tilde{b}^1(\alpha)}
= \alpha (\alpha^2 - 1) P'(\alpha) \frac{\tilde{c}^2(\alpha)}{\tilde{b}^2(\alpha)}
= \dot P^2(\alpha),
\]
where we have used~\eqref{eqn:EMPDi} and~\eqref{eqn:Q reduced} evaluated at $\alpha$.
At the $d_1$ roots of $F^1$, we see that every solution $\dot{P}^2$ takes the same value. Let $\beta$ be such a root, then
\[
\dot{P}^2(\beta)
= \mathbf{\dot{P}}^2(\beta) + 2s(\beta)F^1(\beta)\tilde{P}(\beta)
= \mathbf{\dot{P}}^2(\beta)
= \beta (\beta^2-1) P'(\beta) \frac{\tilde{c}^2(\beta)}{\tilde{b}^2(\beta)},
\]
% where we can be sure that $\tilde{b}^2(\beta) \neq 0$ because it cannot be a root of $b^2$ (if it were, $F\neq 1$). 
However, the other solutions have different values at $\beta$, and this provides the following constraint on the choice of $r$:
\[
\mathbf{\dot{P}}^1(\beta) + 2r(\beta)F^2(\beta)\tilde{P}(\beta) 
% = \mathbf{\dot{P}}^2(\beta) 
= \beta (\beta^2-1) P'(\beta) \frac{\tilde{c}^2(\beta)}{\tilde{b}^2(\beta)}.
\]
This constraint is nontrivial because $\beta$ is not a root of $\tilde{P}$ or $F^2$ by the assumption of the nonsingularity of the spectral curve. As $F^1$ has $d_1$ distinct roots, there are $d_1$ constraints.

Likewise, at the $d_2$ roots of $F^2$, we acquire constraints on the choice of $s$. It is always possible to meet these constraints (because, for example, the degree of $s$ is $d_2$ and there are only $d_2$ roots of $F^2$), so we see that there is a common solution $(\mathbf{\dot{P}}, \mathbf{\dot{b}}^1, \mathbf{\dot{b}}^2)$ to~\eqref{eqn:EMPDi}.

This solution is still not unique; there remains one real parameter. For any real number $s$, we have solutions to~\eqref{eqn:EMPDi} of the form
\[
\dot P = \mathbf{\dot{P}} + 2sP, \;\;\;
\dot b^i = \mathbf{\dot{b}}^i + sb^i
% \labelthis{eqn:P soln}\\
\]
However, this freedom rescales $P$. We have chosen a preferred scaling of $P$, so our choice of $s$ is determined.
% Explicitly, if we were to allow other scalings, the formula for $P$ would be, (cf.~\eqref{eqn:def P})
% \[
% P(t) = r(t) \prod_k (\zeta-\alpha_k(t))(1- \bar{\alpha}_k(t)\zeta),
% \]
% where $\alpha_k$ are the roots inside the unit circle and $r(t)$ is some real function. Then from any solution $\dot{P}$ we can determine the derivatives of the roots $\alpha_k$ at $t=0$. Simply differentiate $P$ and evaluate at $\alpha_k$
% \begin{align*}
% \dot{P} &= \dot{r} \prod_k (\zeta-\alpha_k)(1- \bar{\alpha}_k\zeta) \\
% &\quad+ r(0) \sum_k (-\dot{\alpha}_k + (\dot{\alpha}_k\bar{\alpha}_k+\alpha_k\dot{\bar{\alpha}}_k)\zeta - \dot{\bar{\alpha}}_k\zeta^2) \prod_{m\neq k} (\zeta-\alpha_m)(1- \bar{\alpha}_m\zeta)
% \labelthis{eqn:rescale P} \\
% \dot{P}(\alpha_k) &= -\dot{\alpha}_k(1+\alpha_k\bar{\alpha}_k) \prod_{m\neq k} (\alpha_k-\alpha_m)(1- \bar{\alpha}_m \alpha_k).
% \end{align*}
% Thus we know the values of $\dot{\alpha}_k$, independent of choice of $s$ in our solution~\eqref{eqn:P soln}, because any two solutions $\dot{P}$ differ by a multiple of $P$, which vanishes at every root $\alpha_k$. Alternatively, if we take the lowest order of~\eqref{eqn:rescale P},
% \[
% \dot{P}_0 = \mathbf{\dot{P}}_0 + 2sP_0 = \dot r P_0 + \sum_k (-\dot{\alpha}_k)\prod_{m\neq k} (-\alpha_m),
% \]
% so we may ensure that $r\equiv 1$ by choosing $s$ so that $\dot r = 0$. In short, if we fix a scaling of the spectral curve, then there is a unique solution to~\eqref{eqn:EMPDi}.

Finally there is another necessary condition that must be satisfied by our solution $(\dot{P},\dot{b}^1,\dot{b}^2)$. We must satisfy~\eqref{eqn:residueTangent}, so that~\eqref{eqn:residue} holds along the path. But this condition is satisfied automatically. Observe
\begin{align*}
\dot{P}_1 b^i_0 + P_1 \dot{b}^i_0 - 2(\dot{P}_0 b^i_1 + P_0 \dot{b}^i_1)
% &= P_1\dot{b}^i_0 - 2\dot{P}_0b^i_1 + 3P_1\hat{c}_0 - P_0\dot{b}^i_1 + 2P_1\dot{b}^i_0 \\
&= 3\bra{ P_1\dot{b}^i_0 - \dot{P}_0b^i_1 + P_1\hat{c}_0} \\
&= \frac{3}{P_0}\bra{ P_0P_1\dot{b}^i_0 - P_0\dot{P}_0b^i_1 + P_1\bra{ \frac{1}{2}\dot{P_0}b^i_0 - P_0\dot{b^i_0} }} \\
&= \frac{3\dot{P}_0}{P_0}\bra{ - P_0b^i_1 + \frac{1}{2}P_1b^i_0 }
= 0
\end{align*}
The substitution in the first line comes from the $\zeta^1$ terms of~\eqref{eqn:Q}, the second line from the constant terms of~\eqref{eqn:Q} and the last line comes from the fact that the quantity in the bracket is exactly the residue at $\zeta=0$, which is zero by the assumption that $(P,b^1,b^2)$ lies in $\mathcal{M}_g$, the space of spectral data.

Hence $(\dot{P},\dot{b}^1,\dot{b}^2)$ is a tangent vector to $\mathcal{M}_g$ at $(P,b^1,b^2)$.
\end{proof}
\end{lem}

\begin{lem}[Case (b)]\label{lem:tangent G}
Take a spectral triple $(P,b^1,b^2)\in\mathcal{M}_g$ associated with a nonconformal harmonic map, with a nonsingular spectral curve given by $\eta^2 = P$ of genus $g$. Suppose that $G = \gcd(b^1,b^2)$ is a non-constant real polynomial that does not divide $P$. 
If $G$ lies in $\mathcal{P}^1_\R$ then for every polynomial $\tilde{Q}\in\mathcal{P}^1_\R$, or if $G$ lies in $\mathcal{P}^2_\R$ then for every pair of real numbers $(\tilde{Q},r)$, there exist unique real polynomials $c^i\in\mathcal{P}^{g+1}_\R$ that factor as per~\eqref{eqn:def c Q tilde} and satisfy~\eqref{eqn:Q reduced}. Further, for each such pair $(c^1,c^2)$, there is a unique vector $(\dot P, \dot b^1, \dot b^2) \in \mathcal{P}_\R^{2g+2}\times\mathcal{P}_\R^{g+3}\times\mathcal{P}_\R^{g+3}$ that satisfies~\eqref{eqn:EMPDi} and~\eqref{eqn:residueTangent}.

\begin{proof}
The proof of this lemma is similar to that of~\ref{lem:tangent generic}. We proceed by first solving~\eqref{eqn:Q reduced} and using the resulting pair $(c^1,c^2)$ as inputs to solve~\eqref{eqn:EMPDi}. Regardless of the degree of $G$, which we recall is denoted $d_G$, we must set $Q = G\tilde{Q}$. Equation~\eqref{eqn:Q reduced} reads
\[
\tilde{b}^1\tilde{c}^2 - \tilde{b}^2\tilde{c}^1 = \tilde{Q}\tilde{P}.
\]
There is a unique solution to this equation $(\mathbf{\tilde{c}}^1,\mathbf{\tilde{c}}^2)$ of degree at most $(g+2-d_1-d_G,g+2-d_2-d_G)$. 
If $G$ is linear, multiplying the unique solution of degree at most $(g+1-d_1,g+1-d_2)$ by $F^1$ and $F^2$ respectively gives the desired pair $(c^1,c^2)$.
But if $G$ is quadratic the space of solutions to~\eqref{eqn:Q reduced} is
\[
\Set { (\mathbf{\tilde{c}}^1 + r\tilde{b}^1, \mathbf{\tilde{c}}^1 + r\tilde{b}^2) }
{ r \in \R }.
\]
Hence in that case, for every $r\in\R$ there is a unique pair $(c^1, c^2)$ that factors as required and solves~\eqref{eqn:Q reduced}. In both cases it was not necessary to have a extraneous condition such as $R(Q)=0$, but the choice of $Q$ was restricted by $Q = G\tilde{Q}$.

Next we must solve~\eqref{eqn:EMPDi}, but the proof in Lemma~\ref{lem:tangent generic} applies essentially without modification. 
% The equation in this case is
% \[
% \dot{P} G \tilde{b}^i - 2 F^j \tilde{P} \dot{b}^i = 2 F^j \tilde{P} (\hat{c}^i - \zeta\hat{c}^{i\prime}) + \zeta(\zeta^2-1)P' \tilde{c}^i.
% \labelthis{eqn:EMPDi nc G}
% \]
It has a solution because $\gcd(F^j\tilde{P},G\tilde{b}^i) = 1$.
Analysis at the roots of $F^1F^2\tilde{P}$ shows that there is a common solution $(\mathbf{\dot{P}}, \mathbf{\dot{b}}^1, \mathbf{\dot{b}}^2)$. Again, a choice of scaling of $P$ forces a unique solution. This solution also satisfies~\eqref{eqn:residueTangent}. Hence it is a tangent vector to $\mathcal{M}_g$ at $(b^1,b^2,P)$.
\end{proof}
\end{lem}

Both lemmata above tell essentially the same story, that there there is a choice of two real parameters each giving a unique tangent vector. Conversely, given any tangent vector $(\dot P, \dot{b}^1, \dot{b}^2)$ to $\mathcal{M}_g$ there is a unique pair of polynomials $(\hat{c}^1,\hat{c}^2)$, as shown in Lemma~\ref{lem:unique cs}, and thus a unique polynomial $Q$ from~\eqref{eqn:Q}. Hence this pairing between parameters and tangent vectors is bijective, and we may identify the tangent space to $\mathcal{M}_g$ with these two real parameters. This suggests that $\mathcal{M}_g \cap \mathcal{U}$ itself is a surface.

\begin{thm}\label{thm:moduli manifold}
The open subset $\mathcal{M}_g \cap \mathcal{U}_{ab}$ of the space of spectral triples $\mathcal{M}_g$ is a two dimensional manifold.

\begin{proof}
Recall Definition~\ref{def:subsets U} of $\mathcal{U}_{ab}$ as the open set whose points correspond to cases (a) or (b) and note this is an open set. At any point $p\in \mathcal{M}_g \cap \mathcal{U}_{ab}$, take a simply connected open neighbourhood $\mathcal{V}\subset \mathcal{U}$. On this neighbourhood, define the map $\Psi : \mathcal{V} \to \R^{4g+9}$ by
\begin{align*}
\Psi(P,b^2,b^2) = \Big(
& \int_{A_1} \Theta^1, \dots, \int_{A_g} \Theta^1, \int_{B_1} \Theta^1, \dots, \int_{B_g} \Theta^1, \\
& \int_{A_1} \Theta^2, \dots, \int_{A_g} \Theta^2, \int_{B_1} \Theta^2, \dots, \int_{B_g} \Theta^2, \\
& \int_{\gamma_+} \Theta^1, \int_{\gamma_-} \Theta^1, \int_{\gamma_+} \Theta^2, \int_{\gamma_-} \Theta^2, \\
& P_1b^1_0 - 2P_0b^1_1,\, P_1b^2_0 - 2P_0b^2_1,\, (P_0)^{-1} \prod_{k}(-\alpha_k)
\Big)
\labelthis{eqn:def Psi}
\end{align*}
where $A_i, B_i$ are the real and imaginary periods of $\Sigma$, $\gamma_+$ and $\gamma_-$ are paths in $\Sigma$ between the points over $\zeta=1$ and $\zeta=-1$, and $\alpha_k$ are the roots of $P$ inside the unit circle. Because $\mathcal{V}$ is simply connected, the choice of paths $\{A_k\}, \{B_k\}, \gamma_+, \gamma_-$ may be made smoothly and consistently. The components of $\Psi$ are the conditions that spectral data must satisfy. In particular, the first $4g$ components of $\Psi$ are the periods of the differentials, the next four components are the integrals in the closing conditions~\ref{P:closing}, followed by the conditions to have no residues~\eqref{eqn:residue condition} and the last component of $\Psi$ is our preferred scaling of the spectral curve.

Hence $\mathcal{M}_g \cap \mathcal{V}$ is contained in the level sets of $\Psi$,
\begin{align*}
\mathcal{M}_g \cap \mathcal{V}
\subset \Psi^{-1}\big( &0,\dots,0,2\pi\iu\Z,\dots,2\pi\iu\Z,0,\dots,0,2\pi\iu\Z,\dots,2\pi\iu\Z, \\
&\quad 2\pi\iu\Z,2\pi\iu\Z,2\pi\iu\Z,2\pi\iu\Z,0,0,1 \big).
\end{align*}
The point $p$ of $\mathcal{M}_g \cap \mathcal{V}$ falls under either Lemma~\ref{lem:tangent generic} or~\ref{lem:tangent G}. In both cases we computed that the kernel of $d\Psi_p$ is two dimensional. The differential of $\Psi$ is a map from $\R^{4g+11}$ to $\R^{4g+9} = \R^{4g+4+2+2+1}$, and so is full rank at every such point $p$. Therefore by the Implicit Function Theorem $\mathcal{M}_g \cap \mathcal{U}_{ab}$ is a two dimensional manifold.
\end{proof}
\end{thm}

%%%%%%%%%%%%%%%%%%%%%%%%%%%%%%%%%%%%%%%%%%%%%%%%%%%%%%%%%%%%%%%%%%%%
%%%%%%%%%%%%%%%%%%%%%%%%%%%%%%%%%%%%%%%%%%%%%%%%%%%%%%%%%%%%%%%%%%%%
%%%%%%%%%%%%%%%%%%%%%%%%%%%%%%%%%%%%%%%%%%%%%%%%%%%%%%%%%%%%%%%%%%%%
%%%%%%%%%%%%%%%%%%%%%%%%%%%%%%%%%%%%%%%%%%%%%%%%%%%%%%%%%%%%%%%%%%%%
\section{Conformal Harmonic Maps}\label{sec:conformal}

From equation~\eqref{eqn:EMPDi all factors} onwards, we made the assumption that the spectral triple came from a nonconformal harmonic map. However, conformal harmonic maps are of particular interest because they are minimal immersions. We show in this section, under mild assumptions, that the points of $\mathcal{M}_g$ corresponding to conformal maps are also smooth points of dimension two. 

Let us return to the discussion following~\eqref{eqn:EMPDi all factors}, this time assuming that a point of $\mathcal{M}_g$ corresponds to a conformal map. We know that $P_0(t)$, $b^1_0(t)$ and $b^2_0(t)$ all vanish at $t=0$. Thus $F=\gcd(P,b^1,b^2)$ includes a factor of $\zeta$.
It is now $\zeta^{-1}FF^i$ that divides $c^i$ and so $\zeta^{-1}FG$ divides $Q$, by~\eqref{eqn:Q}. However the residue condition~\eqref{eqn:residueTangent} in this case simplifies in a way that forces another constraint on $Q$. At $t=0$,~\eqref{eqn:residueTangent} becomes
\[
P_1 \dot{b}_0 - 2 \dot{P}_0 b_1 = 0.
\]
Combining this with the terms of linear degree in equation~\eqref{eqn:EMPDi} gives
% \[
% \dot P_0 b_1^i - 2P_1\dot b_0^i = 3P_1\hat{c}_0^i.
% \]
% Combining these two expressions shows that
\[
3P_1\hat{c}_0^i = \dot P_0 b_1^i - 4\dot{P}_0 b_1^i = -3\dot{P}_0 b_1^i.
\]
Substituting this into the linear degree of~\eqref{eqn:Q}, we arrive at
\begin{align*}
% Q_0 P_1 &= b^1_1 c^2_0 - b^2_1 c^1_0 \\
Q_0 (P_1)^2 
&= b^1_1 (P_1 c^2_0) - b^2_1 (P_1 c^1_0)
% &= b^1_1 (\dot{P}_0 b^2_1) - b^2_1 (\dot{P}_0 b^1_1) 
= 0.
\labelthis{eqn:Q0 vanish conformal}
\end{align*}
A spectral curve must be nonsingular at $\zeta=0$, so if $P_0=0$ we can be sure that $P_1\neq 0$. Hence $Q_0$ must vanish. As $Q$ is a real quadratic polynomial, it must be of the form $Q=Q_1 \zeta$ for some real number $Q_1$. Immediately it follows that if a deformation exists at a point corresponding to a conformal map then $F = \zeta$ and $G = 1$, as the polynomials $b^i$ are not permitted to have multiple roots at $\zeta=0$. Thus there are two conformal cases, and only (e) can lead to deformations.

\begin{center}
\begin{tabular}{|c|c|c|c|}
\hline
Case & $P_0$ & $\deg F$ & $\deg G$ \\ \hline\hline
(e) & \multirow{2}{*}{$P_0 = 0$} & $F=\zeta$ & 0 \\ \cline{1-1}\cline{3-4}
(f) && \multicolumn{2}{|c|}{$FG \in \mathcal{P}^k_\R,\, k>2$} \\ \hline
\end{tabular}
\end{center}

Without further delay, let us show that the necessary equations can be solved in case (e).

\begin{lem}[Case (e)]\label{lem:tangent conformal}
Take a spectral triple $(P,b^1,b^2)\in\mathcal{M}_g$ associated with a conformal harmonic map, with a nonsingular spectral curve given by $\eta^2 = P$ of genus $g$. 
Suppose that $\gcd(b^1,b^2) = \zeta$. 
Then for every pair of real numbers $(Q_1,r)$, there exist unique real polynomials $c^i\in\mathcal{P}^{g+1}_\R$ that factor as per~\eqref{eqn:def c Q tilde} and satisfy~\eqref{eqn:Q reduced}. Further, for each such pair $(c^1,c^2)$, there is a unique vector $(\dot P, \dot b^1, \dot b^2) \in \mathcal{P}_\R^{2g+2}\times\mathcal{P}_\R^{g+3}\times\mathcal{P}_\R^{g+3}$ that satisfies~\eqref{eqn:EMPDi} and~\eqref{eqn:residueTangent}. It is therefore a tangent vector to the space of spectral data $\mathcal{M}_g$.

\begin{proof}
This is the conformal case, so $P(0) = P_0 = 0$. From~\eqref{eqn:residue condition}, $b^i_0 = 0$ also. We may write therefore that $P= \zeta F^1F^2\tilde{P}$ and $b^i = \zeta F^i \tilde{b}^i$, where $\tilde{P}\in\mathcal{P}^{2g-d_1-d_2}$ and the polynomials $\tilde{b}^i\in\mathcal{P}^{g+1-d_i}_\R$ are coprime.

We have already demonstrated in~\eqref{eqn:Q0 vanish conformal} that $\zeta$ necessarily divides $Q$. Thus~\eqref{eqn:Q reduced} is simply
\[
\tilde{b}^1 \tilde{c}^2 - \tilde{b}^2 \tilde{c}^1 = \zeta Q_1\tilde{P}.
\labelthis{eqn:Q conformal}
\]
This is similar to the above case where $G$ was quadratic (Lemma~\ref{lem:tangent G}). The space of solutions is
\[
\Set { (\mathbf{\tilde{c}}^1 + r \tilde{b}^1, \mathbf{\tilde{c}}^2 + r \tilde{b}^2) }{ r\in \R },
\]
where $(\mathbf{\tilde{c}}^1,\mathbf{\tilde{c}}^2)$ is the unique solution of degree at most $(g+1-d_1, g+1-d_2)$. For every such solution, let $c^i = F^i\tilde{c}^i$ and consider the corresponding~\eqref{eqn:EMPDi}.
% , namely
% \[
% \dot{P} \tilde{b}^i - 2 F^j \tilde{P} \dot{b}^i = 2 F^j \tilde{P} (\hat{c}^i - \zeta\hat{c}^{i\prime}) + (\zeta^2-1)P' \tilde{c}^i.
% \labelthis{eqn:EMPDi conformal}
% \]
As before, there is a common solution $(\mathbf{\dot{P}}, \mathbf{\dot{b}}^1, \mathbf{\dot{b}}^2)$ for the dotted quantities. The space of solutions is however
\[
\Set{
(\mathbf{\dot P} + 2s F^1F^2\tilde{P}, \mathbf{\dot b}^1 + s F^1\tilde{b}^1, \mathbf{\dot b}^2 + s F^2\tilde{b}^2)
}{ s \in\mathcal{P}^2_\R },
\]
and to each choice $(Q,r)$ there are many tangent vectors. However, unlike the cases (a) and (b), equation~\eqref{eqn:residueTangent} is not automatically satisfied. Let $s = s_0 + s_1\zeta + \bar{s}_0 \zeta^2$. For $i=1$, we see that the condition implies that
\[
2 \mathbf{\dot{P}}_0 b^1_1 - P_1 \mathbf{\dot{b}}^1_0 + 3 s_0 P_1 b^1_1 = 0,
\]
which fully determines $s_0$. We now show that this solution simultaneously satisfies the condition for $i=2$. Note that~\eqref{eqn:EMPDi} in the lowest degree reads
\[
\dot{P}_0 b^i_1 - 2P_1\dot{b}^i_0 = -3 P_1 c^i_0,
\]
and~\eqref{eqn:Q conformal} in the lowest degree yields
\begin{align*}
% b^1_1 c^2_0 &= b^2_1 c^1_0 \\
% b^1_1 \bra{\mathbf{\dot{P}_0} b^2_1 - 2P_1\mathbf{\dot{b}^2_0}} &= b^2_1 \bra{\mathbf{\dot{P}_0} b^1_1 - 2P_1\mathbf{\dot{b}^1_0}} \\
2b^1_1 \mathbf{\dot{b}^2_0} &= 2 b^2_1 \mathbf{\dot{b}^1_0}.
\end{align*}
Condition~\eqref{eqn:residueTangent} for $i=2$ is therefore
\begin{align*}
b^1_1 \bra{2 \mathbf{\dot{P}_0} b^2_1 - P_1 \mathbf{\dot{b}^2_0} + 3 s_0 P_1 b^2_1}
% &= 2 \mathbf{\dot{P}_0} b^1_1b^2_1 - P_1 b^1_1\mathbf{\dot{b}^2_0} + 3 s_0 P_1 b^1_1b^2_1 \\
% &= 2 \mathbf{\dot{P}_0} b^1_1b^2_1 + P_1 b^2_1 \mathbf{\dot{b}^1_0} + 3 s_0 P_1 b^1_1b^2_1 \\
&= b^2_1\bra{2 \mathbf{\dot{P}_0} b^1_1 - P_1 \mathbf{\dot{b}^1_0} + 3 s_0 P_1 b^1_1 } = 0.
\end{align*}
Hence we have demonstrated that the condition holds for $i=2$ also. Having cleared this hurdle, there is still one free parameter. For any $Q_1$ and $r$, the corresponding tangent vectors that solve~\eqref{eqn:EMPDi} are
\begin{align*}
\Big\{
\Big(&\mathbf{\dot P} + 2(s_0+\bar{s}_0\zeta^2) F^1F^2\tilde{P} + 2s_1\zeta F^1F^2\tilde{P}, \\
&\qquad\mathbf{\dot b^1} + (s_0+\bar{s}_0\zeta^2) F^1\tilde{b}^1 + s_1\zeta F^1\tilde{b}^1, \\
&\qquad\mathbf{\dot b^2} + (s_0+\bar{s}_0\zeta^2) F^2\tilde{b}^2 + s_1\zeta F^2\tilde{b}^2 \,\Big)
\qquad\Big\vert s_1 \in \R
\Big\}.
\end{align*}
But our free choice of $s_1\in\R$ is only adding multiples of $(2P,b^1,b^2)$, which as in the nonconformal case is a rescaling of the spectral curve, and so also determined uniquely.
\end{proof}
\end{lem}

From this Lemma there is a result analogous to Theorem~\ref{thm:moduli manifold}. As in that theorem, we can show that every point of $\mathcal{M}_g \cap\, \mathcal{U}_{e}$ is a smooth point of $\mathcal{M}_g$. However, $\mathcal{U}_{e}$ is not an open set, so we cannot say that $\mathcal{M}_g \cap\, \mathcal{U}_{e}$ is a manifold. 

\begin{thm}\label{thm:conformal moduli manifold}
The points of $\mathcal{M}_g \cap \mathcal{U}_{e}$, corresponding to conformal harmonic maps, are smooth points of $\mathcal{M}_g$ of dimension two.

\begin{proof}
Is entirely similar to Theorem~\ref{thm:moduli manifold}. 
% At such a point $p\in\mathcal{M}_g \cap \mathcal{U}_{e}$ we have that $P(0) = 0$ and $\gcd(b^1,b^2)=\zeta$. Take a simply connected open neighbourhood $\mathcal{V}\subset \mathcal{U}$ of $p$ as before. We may use the same definition of $\Psi$, equation~\eqref{eqn:def Psi}, and by the implicit function theorem and Lemma~\ref{lem:tangent conformal} it follows that $p$ is also a smooth point of $\mathcal{M}_g$.
\end{proof}
\end{thm}

% A final point to note about conformal spectral triples within $\mathcal{M}_g$. Suppose there was a deformation that stayed within the space of conformal spectral triples. Then $\dot{P}_0$ and $\dot{b}^i_0$ would both be zero, so~\eqref{eqn:} would be automatically satisfied, and $\zeta$ would divide $\dot{P}$ and $\dot{b}^i$. From~\eqref{eqn:EMPDi conformal}, this forces $\zeta$ to divide $\tilde{c}^i$ also.

%%%%%%%%%%%%%%%%%%%%%%%%%%%%%%%%%%%%%%%%%%%%%%%%%%%%%%%%%%%%%%%%%%%%
%%%%%%%%%%%%%%%%%%%%%%%%%%%%%%%%%%%%%%%%%%%%%%%%%%%%%%%%%%%%%%%%%%%%
%%%%%%%%%%%%%%%%%%%%%%%%%%%%%%%%%%%%%%%%%%%%%%%%%%%%%%%%%%%%%%%%%%%%
%%%%%%%%%%%%%%%%%%%%%%%%%%%%%%%%%%%%%%%%%%%%%%%%%%%%%%%%%%%%%%%%%%%%
\section{Consequences}

In this final section we collect some corollaries. First we show that when the genus $g$ of the spectral curve is zero or one, it cannot have singularities and cases (c), (d) and (f) do not occur. Thus $\mathcal{M}_0$ and $\mathcal{M}_1$ are smooth at every point and therefore are surfaces. Then we comment on why case (c) is difficult. And finally we give an interpretation of $Q$.

Firstly, we can rule out singular spectral curves with genus two or less. In (arithmetic) genus zero, there is one one pair of branch points so singularities are plainly impossible. In higher genus, suppose we have a singular spectral curve $\Sigma$ with normalisation $\tilde{\Sigma}$.
Because there can be no singular points on the unit circle, all singular points come in pairs, so the genus of $\Sigma$ and $\tilde{\Sigma}$ differ by at least two. This excludes the possibility of singular spectral curves of genus zero or one. If we add in the fact that at a point of $\tilde{\Sigma}$ that maps to $\Sigma$ with multiplicity $m$, the differentials $\Theta$ and $\tilde{\Theta}$ both have a common zero of order $m-1$, then Lemma~\ref{lem:no singularities} below shows that no singular spectral curve has a genus zero normalisation. This rules out singular genus two spectral curves also.

To exclude the possibility of the undesirable cases (c), (d) and (f), we need to know if the differentials have a common root. For genus zero this is answered by the following lemma. 
In~\cite{Carberry2016}, a stronger result analogous to~\ref{lem:no singularities} is proved by examining the degree of $b^1/b^2$ considered as a function $\CP^1 \to \CP^1$.

\begin{lem}\label{lem:no singularities}
On a spectral curve of genus zero, differentials satisfying conditions~\ref{P:poles}--\ref{P:reality} with linearly independent principal parts do not have common roots.
\begin{proof}
We distinguish between two cases: whether or not the spectral curve is branched over $\zeta=0$. If $\zeta=0$ is a branch point, then we note the following more general proof. Suppose the spectral curve has genus $g$. Then from~\eqref{eqn:def b} we have that any differential may be written as 
\[
a(\zeta)\frac{d\zeta}{\zeta\eta},
\]
for some real polynomial $a$ of degree $g+1$. Any real polynomial is determined up to real scaling by its $g+1$ roots, so if two such differentials have $g+1$ roots in common, then they are linearly dependent over $\R$. Letting $g=0$ shows that two differentials may not share any roots.

In the nonconformal case, instead we consider the specific form of the differentials. Suppose that the spectral curve is branched over $\alpha$ and $\cji{\alpha}$ and let $x = -\frac{1}{2}\alpha^{-1}(1+\alpha\bar{\alpha})$. Then a differential $\Theta$ satisfying~\ref{P:poles}--\ref{P:reality} is given by
\[
\Theta = (y + xy \zeta + \overline{xy}\zeta^2 + \bar{y}\zeta^3)\frac{d\zeta}{\zeta^2\eta},
\]
for some a nonzero constant $y$. If $\Theta$ has a root at $\beta$, then 
\[
\frac{y}{\bar{y}} = -\beta^2 \frac{\bar{x} + \beta}{1+x\beta}.
\]
So any two such differentials with a common root are linearly dependent over $\R$.
\end{proof}
\end{lem}

We shall now give a similar short proof to show that the differentials of a genus one spectral curve may not have a common root at the branch points, which excludes (c). Suppose that $\Sigma$ is genus one with branch points at $\alpha,\beta,\cji{\alpha},\cji{\beta}$, none of which are zero, and that $\gcd(P,b^1,b^2) = F$ is quadratic. Without loss of generality, let $F = (\zeta-\alpha)(1-\bar{\alpha}\zeta)$ and $b^1 = (c + d\zeta + \bar{c}\zeta^2)F$ for some complex number $c$ and real number $d$. Expanding this and applying~\eqref{eqn:residue condition} shows that
\begin{align*}
0 
% &= \bra{ -\alpha(1+\beta\bar{\beta}) - \beta(1+\alpha\bar{\alpha}) }(-\alpha c) - 2(\alpha\beta)\bra{-\alpha d + (1+\alpha\bar{\alpha})c}\\
&= \bra{ \alpha^2(1+\beta\bar{\beta}) - \alpha\beta(1+\alpha\bar{\alpha}) }c + 2\alpha^2\beta d.
\end{align*}
The coefficient of $c$ above, the bracketed expression, is never zero and so $c$ is determined by $d$. Hence $b^1$ is determined up to a real scalar. This demonstrates any two differentials with the same factor $F$ are real linearly dependent, which contradicts~\ref{P:linear independence}. If the differentials have common roots as per case (d) or (f), then they are forced to be real linearly dependent by the residue condition.

Thus we have shown that for a spectral curve of genus zero or one, it cannot have singularities and cases (c), (d) and (f) do not occur. Thus $\mathcal{M}_0$ and $\mathcal{M}_1$ are smooth at every point and therefore are surfaces.

In higher genus, it is possible that there are exist harmonic tori that fall under cases (c), (d), or (f). Cases (d) and (f) do not admit deformations, so could be isolated points of $\mathcal{M}_g$, cusps, or some other type of singular point. Case (c) however could be a smooth point of $\mathcal{M}_g$. The obstacle to proving a result analogous to Lemma~\ref{lem:tangent G} is solving equation~\eqref{eqn:Q reduced}, which would read
\[
\tilde{b}^2\tilde{c}^2 - \tilde{b}^2\tilde{c}^1 =\tilde{Q}\tilde{P},
\]
for $\tilde{Q}\in\R$. As the degree of $\tilde{b}^2$ is $g+1-d_2$, there is a solution such that $\tilde{c}^2$ has degree $g-d_2$. However, if $\tilde{c}^2$ came from a deformation, then we would expect it to have degree $g-1-d_2$. This is similar to case (a), where we imposed a condition on $Q$ to ensure the correct degree, namely that $R(Q)=0$. In this case however, the space of allowed $Q$ is already one dimensional. Therefore whether or not an infinitesimal deformation exists depends in a transcendental way on the spectral curve.

We can spell out this dependence explicitly. Given a spectral curve, it determines a plane of differential satisfying conditions~\ref{P:poles}--\ref{P:periods} and therefore determines $F=\gcd(P,b^1,b^2)$. The space $\mathcal{M}_g\cap \mathcal{U}_{e}$ is the subset of spectral triples where $F$ is quadratic and the plane of differentials contains two linearly independent differentials with integral periods. To show that deformations are always possible in case (c) is equivalent to showing that $R((P,b^1,b^2),F)$ vanishes on $\mathcal{M}_g\cap \mathcal{U}_{e}$.

To close, let us give a geometric interpretation to the polynomial $Q$. The conformal type of the domain of a harmonic tori is given by the ratios of the principal parts of the differentials of its spectral data. Let the conformal type be denoted $\tau$. For a nonconformal harmonic tori we have that $b^2_0 = \tau b^1_0$. Consideration of the constant terms of~\eqref{eqn:EMPDi} reveal that
\[
\dot{P}_0b^i_0 -2P_0 \dot{b}^i = 2P_0\hat{c}^i_0.
\]
% Noting that $c^i_0 = - \hat{c}^i_0$, we
We substitute this into~\eqref{eqn:Q} to arrive at
\[
Q_0 P_0 
% = b^1_0 c^2_0 - b^2_0 c^1_0 
= b^1_0 \dot{b}^2_0 - \dot{b}^1_0 b^2_0.
\]
Differentiating the relationship $b^2_0 = \tau b^1_0$ and rearranging the above yields
\[
Q_0 = \frac{\dot{\tau}}{\tau} \frac{b^1_0 b^2_0}{P_0}.
\labelthis{eqn:Q0 change conformal}
\]
We see therefore that $Q_0$ controls the change in the conformal type of the domain of the harmonic map.

\bibliographystyle{ross_bibsty}
\bibliography{zotero}

\newcommand{\etalchar}[1]{$^{#1}$}
\begin{thebibliography}{GKM{\etalchar{+}}95}

\bibitem[CS16a]{Carberry2016a}
Carberry E and Schmidt M~U.
\newblock The closure of spectral data for constant mean curvature tori in
  𝕊3.
\newblock \emph{Journal f{\"u}r die reine und angewandte Mathematik (Crelles
  Journal)}, 2016(721), 2016.

\bibitem[CS16b]{Carberry2016}
Carberry E and Schmidt M~U.
\newblock The prevalence of tori amongst constant mean curvature planes in r3.
\newblock \emph{Journal of Geometry and Physics}, 106:352--366, 2016.

\bibitem[FFM80]{Flaschka1980}
Flaschka H, Forest M~G, and McLaughlin D~W.
\newblock Multiphase averaging and the inverse spectral solution of the
  {{Korteweg}}-de {{Vries}} equation.
\newblock \emph{Communications on Pure and Applied Mathematics},
  33(6):739--784, 1980.

\bibitem[GKM{\etalchar{+}}95]{Gorsky1995}
Gorsky A, Krichever I, Marshakov A, Mironov A, and Morozov A.
\newblock Integrability and {{Seiberg}}-{{Witten}} exact solution.
\newblock \emph{Physics Letters B}, 355(3-4):466--474, 1995.

\bibitem[Gru12]{Grushevsky2012}
Grushevsky S.
\newblock The {{Schottky}} problem.
\newblock In \emph{Current Developments in Algebraic Geometry}, volume~59 of
  \emph{Mathematical Sciences Research Institute Publications}, pages 129--164.
  {Cambridge University Press}, Cambridge, 2012.

\bibitem[Hit90]{Hitchin1990}
Hitchin N~J.
\newblock Harmonic maps from a 2-torus to the 3-sphere.
\newblock \emph{Journal of Differential Geometry}, 31:627--710, 1990.

\bibitem[HKS16]{Schmidt2016}
Hauswirth L, Kilian M, and Schmidt M~U.
\newblock Singularities of {{Whitham}} flows for hyperelliptic spectral curves,
  2016.
\newblock Preprint.

\bibitem[Kal84]{Kalman1984}
Kalman D.
\newblock The generalized {{Vandermonde}} matrix.
\newblock \emph{Mathematics Magazine}, 57(1):15--21, 1984.

\bibitem[Kri95]{Krichever1995}
Krichever I.
\newblock Algebraic-geometrical methods in the theory of integrable equations
  and their perturbations.
\newblock \emph{Acta Applicandae Mathematicae}, 39(1-3):93--125, 1995.

\bibitem[KSS15]{Kilian2015a}
Kilian M, Schmidt M~U, and Schmitt N.
\newblock Flows of constant mean curvature tori in the 3-sphere: The
  equivariant case.
\newblock \emph{Journal f{\"u}r die Reine und Angewandte Mathematik (Crelle's
  Journal)}, 707:45--86, 2015.

\bibitem[LL83]{Lax1983}
Lax P~D and Levermore C~D.
\newblock The small dispersion limit of the {{Korteweg}}-de {{Vries}} equation.
  {{I}}.
\newblock \emph{Communications on Pure and Applied Mathematics},
  36(10):253--290, 1983.

\bibitem[Mir95]{Miranda1995}
Miranda R.
\newblock \emph{Algebraic Curves and {{Riemann}} Surfaces}, volume~5 of
  \emph{Graduate Studies in Mathematics}.
\newblock {American Mathematical Society, Providence, RI}, 1995.

\bibitem[Shi86]{Shiota1986}
Shiota T.
\newblock Characterization of {{Jacobian}} varieties in terms of soliton
  equations.
\newblock \emph{Inventiones Mathematicae}, 83(2):333--382, 1986.

\bibitem[Whi72]{Whitney1972}
Whitney H.
\newblock \emph{Complex Analytic Varieties}.
\newblock {Addison-Wesley Publishing Co.}, Reading, Massachusetts, 1972.

\end{thebibliography}
\end{document}